\documentclass[oneside]{amsart}

\usepackage[utf8]{inputenc}
\usepackage[T1]{fontenc}

\usepackage[
  left=2.5cm,
  right=2.5cm,
  top=2.5cm,
  bottom=2.5cm
]{geometry}

\usepackage{
  amsmath, amsthm, amssymb,
  mathrsfs,
  fancyhdr,
  graphicx,
  url,
  caption,
  float,
  scrextend,
  array,
  pifont,
  soul,
  textcomp,
  xcolor,
  wrapfig,
  enumitem,
  stackengine,
  mathtools,
  tikz, tikz-cd
}

\usepackage{tikz, tikz-cd}
\usepackage[
  colorlinks=true,
  urlcolor=black,
  citecolor=black,
  linkcolor=black,
  hyperfootnotes=true
]{hyperref}

\stackMath
\newcommand\tsup[2][2]{%
 \def\useanchorwidth{T}%
  \ifnum#1>1%
    \stackon[-1pt]{\tsup[\numexpr#1-1\relax]{#2}}{\hspace{1pt}\scriptstyle\sim}%
  \else%
    \stackon[.5pt]{#2}{\hspace{1pt}\scriptstyle\sim}%
  \fi%
}

\newcommand{\nc}{\newcommand}
\nc{\Opn}{\mathrm{O}}

\nc{\cO}{\mathcal{O}}

\newcommand{\sgo}{{\mathsf{S}_1(\Ga,\Op)}}

\nc{\ram}{\mathsf{Ramsey}}
\nc{\soo}[1]{\mathsf{S}_1(\Op,\Op)}
\nc{\swl}[1]{\mathsf{S}_1(\Om,\Lambda)}
\nc{\swg}[1]{\mathsf{S}_1(\Om,\Ga)}
\nc{\goo}[1]{\gone(\Op,\Op)}
\nc{\gwo}{\gone(\Om,\cO)}
\nc{\gwl}[1]{\gone(\Om(#1),\Lambda(#1))}
\nc{\soox}[1]{\mathsf{S}_1(\cO(#1),\cO(#1))}
\nc{\swlx}[1]{\mathsf{S}_1(\Om(#1),\Lambda(#1))}
\nc{\goox}[1]{\gone(\cO(#1),\cO(#1))}
\nc{\gwox}[1]{\gone(\Om(#1),\cO(#1))}
\nc{\gwlx}[1]{\gone(\Om(#1),\Lambda(#1))}
\nc{\mich}[1]{{(#1)_\mathrm{M}}}
\nc{\michk}[2]{{(#1)^{#2}_\mathrm{M}}}
\nc{\PNM}{{\PN_\mathrm{M}}}
\nc{\sfinwl}{\mathsf{S}_{\mathrm{fin}}(\Om,\Lambda)}
\nc{\gfinoo}{\gfin(\cO,\cO)}
\nc{\gfinwo}{\gfin(\Om,\cO)}
\nc{\gfinwl}{\gfin(\Om,\Lambda)}
\nc{\sfinoox}[1]{\mathsf{S}_{\mathrm{fin}}(\cO(#1),\cO(#1))}
\nc{\sfinwlx}[1]{\mathsf{S}_{\mathrm{fin}}(\Om,\Lambda)}
\nc{\sfinwwx}[1]{\mathsf{S}_{\mathrm{fin}}(\Om(#1),\Om(#1))}
\nc{\gfinoox}[1]{\gfin(\cO(#1),\cO(#1))}
\nc{\gfinwox}[1]{\gfin(\Om(#1),\cO(#1))}
\nc{\gfinwlx}[1]{\gfin(\Om(#1),\Lambda(#1))}

\nc{\mc}{\mathcal}
\nc{\thusfar}{\my{--- Edited thus far ---}}
\nc{\lei}{\le^\oo}
\nc{\sqsubs}{\sqsubseteq^*}
\nc{\card}[1]{\left|#1\right|}
\nc{\medcard}[1]{\biggl|\,#1\,\biggr|}
\nc{\smallmedcard}[1]{\bigl|\,#1\,\bigr|}
\nc{\smallcard}[1]{|\,#1\,|}
\nc{\bds}{bidirectional $\roth$-scale}
\nc{\bfP}{\mathbf{P}}
\nc{\bfQ}{\mathbf{Q}}
\nc{\bbT}{\mathbb{T}}
\nc{\bbZ}{\mathbb{Z}}
\nc{\bbN}{\mathbb{N}}
\nc{\bbC}{\mathbb{C}}
\nc{\beq}{\begin{equation}}
\nc{\eeq}{\end{equation}}
\nc{\beqs}{\begin{equation*}}
\nc{\eeqs}{\end{equation*}}
\nc{\mbq}{\mb{?}}
\nc{\mb}[1]{{\mbox{\textbf{#1}}}}
\nc{\nop}{$\times$}
\nc{\fbn}{\!\!\fbox{\!\nop\!}\!\!}
\nc{\yup}{\checkmark}
\nc{\forces}{\Vdash}
\nc{\name}[1]{\dot{#1}}
\nc{\tf}{\my{FINISHED THUS FAR}}
\nc{\FU}{Fr\'echet--Urysohn}
\nc{\gs}{$\gamma$~space}
\nc{\Gab}{\Gamma_{\mathrm{B}}}
\nc{\Omb}{\Omega_{\mathrm{B}}}
\nc{\Ga}{\Gamma}
\nc{\Om}{\Omega}
\nc{\smallbinom}[2]{\begin{psmallmatrix} #1\\ #2 \end{psmallmatrix}}
\nc{\bgamma}{\smallbinom{\Om}{\Ga}}

\nc{\productive}[2]{(#1,\allowbreak #2)^\x}
\nc{\prdct}[1]{#1^\x}
\nc{\Sel}{\mathsf{S}}
\nc{\sset}[2]{\{\,#1 : #2\,\}}
\nc{\smb}[1]{{\!\!\mb{#1}\!\!}}
\nc{\medset}[2]{{\biggl\{\,#1 : #2\,\biggr\}}}
\nc{\smallmedset}[2]{{\bigl\{\,#1 : #2\,\bigr\}}}
\nc{\set}[2]{{\left\{\,#1 : #2\,\right\}}}
\nc{\seq}[2]{{\la\, #1 : #2\,\ra}}
\nc{\eseq}[1]{#1_0, \allowbreak #1_1, \allowbreak\dotsc} 

\nc{\eseqint}[3]{#1_{#2}, \allowbreak\dotsc,\allowbreak #1_{#3}} 
\nc{\eseqstart}[2]{#1_{#2},\allowbreak #1_{#2+1},\dotsc } 
\nc{\eprod}[1]{#1_{1}\times \allowbreak#1_{2}\times\dotsb}

\nc{\shortprod}[1]{\prod_{n=1}^\infty{#1}_n}

\nc{\eprodint}[3]{#1_{#2}\times \allowbreak\dotsb\times\allowbreak #1_{#3}}

\nc{\seleseq}[1]{#1_1\in \mathcal{#1}_1, \allowbreak #1_2\in \mathcal{#1}_2, \allowbreak\dotsc}
\nc{\cube}{(\Cantor)^\bbN}
\nc{\Match}{\op{Match}}
\nc{\concat}[1]{\hat{\phantom{a}}\langle #1\rangle}
\nc{\poset}{\mathbb{P}}
\nc{\fn}[1]{{\op{Fn}(#1\times\w,2)}}
\nc{\linadd}{\op{linadd}}
\nc{\nonprod}{\non^\x}
\nc{\alephes}{{\aleph_0}}
\nc{\my}[1]{{\color{red}{#1}}}
\nc{\later}[1]{{\color{green} #1}}
\nc{\BTs}[1]{{\color{green} #1 (BT)}}
\nc{\Cp}{\op{C}_\mathrm{p}}
\nc{\Bp}{\op{B}_p}
\nc{\Pa}[8]{\bibitem{#1} {#2}, \emph{#3}, {#4} \textbf{#5} ({#6}), {#7}--{#8}.}
\nc{\tPa}[5]{\bibitem{#1} {#2}, \emph{#3}, {#4}, to appear.}
\nc{\sPa}[4]{\bibitem{#1} {#2}, \emph{#3}, {#4}, submitted.}
\nc{\Bc}[9]{\bibitem{#1} {#2}, \emph{#3}, in: \textbf{#4} (#5), #6 #7, #8--#9.}
\nc{\fD}{\mathfrak{D}}
\nc{\fX}{\mathfrak{X}}
\nc{\Onbd}{\Op_{\mathrm{nbd}}} 
\nc{\Omnb}{\Om_{\mathrm{nbd}}} 
\nc{\od}{\mathfrak{od}}
\nc{\Setting}[7]{\xymatrix@R=4pt@C=7pt{#1\ar@{-}[r]&#2\ar@{-}[r]&#3\\&#4\ar@{-}[u]\\
#5\ar@{-}[uu]\ar@{-}[r] & #6\ar@{-}[u]\ar@{-}[r] & #7\ar@{-}[uu]}}
\nc{\mx}[1]{\begin{matrix}#1\end{matrix}}
\nc{\plim}{p\txt{-}\lim}
\nc{\Bgp}{{\Z^\bbN}}
\nc{\Cgp}{{{\Z_2}^\bbN}}
\nc{\Cite}[1]{\textbf{[#1]}}
\nc{\Next}[1]{{#1^+}}
\nc{\cFin}{\mathrm{cF}}
\nc{\scsp}{\text{-scale space}}
\nc{\cfn}{\text{cofinal}\ }
\nc{\Con}{\text{Concentrated}}
\nc{\Lind}{\text{Lindel\"of}\,}
\nc{\con}{\text{-Concentrated}}
\nc{\lind}{\text{-Lindel\"of}\,}
\nc{\ctbl}{\text{countably }\allowbreak}
\nc{\Hur}{\text{Hurewicz}}
\nc{\intvl}[2]{{[#1(#2),\allowbreak #1(#2\!+\!1))}}
\nc{\Bdd}{\mathbf{B}}
\nc{\Dfin}{\mathfrak{D}_\mathrm{fin}}
\nc{\grbl}{{\mbox{\textit{\tiny gp}}}}
\nc{\bbP}{\mathbb{P}}
\nc{\BOfat}{\B_{\Om_{\mathrm{fat}}}}
\nc{\Bgood}{\B_{\mathrm{good}}}
\nc{\compactN}{\cl{\mathbb{N}}}
\nc{\blocks}[2]{\op{cl}_{#2}(#1)}
\nc{\blocksplus}[2]{\op{cl}^+_{#2}(#1)}
\nc{\arx}[1]{\texttt{http://arxiv.org/math/#1}}
\nc{\bq}{\begin{quote}}
\nc{\eq}{\end{quote}}
\nc{\cl}[1]{\overline{#1}}
\nc{\Cl}[2]{\mathrm{cl}_{#1}(#2)}
\nc{\CH}{the Continuum Hypothesis}
\nc{\MA}{Martin's Axiom}
\nc{\Bfat}{\B_\mathrm{fat}}
\nc{\inv}{^{-1}}
\nc{\Cantor}{{2^\w}}
\nc{\bP}{\mathbf{P}}
\nc{\bof}{\op{\fb}}
\nc{\dof}{\op{\fd}}
\nc{\bofF}{\bof(\cF)}
\nc{\sr}[3]{\underset{\mbox{#3}}{\mbox{#1}}}
\nc{\gp}{\binom{\Om}{\Ga}}
\nc{\gpsmall}{\mbox{$\gp$}}
\nc{\gig}{\gimel}
\nc{\gns}{\sone(\Om,\gig)}
\nc{\nsr}[2]{#1}
\nc{\Srg}{{\mathbb{S}}}
\nc{\Srgs}{{\mathbb{S}^*}}
\nc{\NN}{{\w^{\w}}}
\nc{\ZN}{{\Z^{\bbN}}}
\nc{\NNup}{{\bbN^{\uparrow\bbN}}}
\nc{\NNupb}{{b^{\uparrow\bbN}}}
\nc{\Pof}{\op{P}}
\nc{\PN}{{\Pof(\w)}}
\nc{\rothx}[1]{{[#1]^{\mbox{\tiny $\infty$}}}}
\nc{\tx}{{\tilde{x}}}
\nc{\roth}{{[\w]^{\w}}}
\nc{\roths}{{[b]^{\mbox{\tiny $\infty$}}}} 
\nc{\Fin}{\mathrm{Fin}}
\nc{\ici}{[\w]^{ \w, \w}}
\nc{\icikap}{[\kappa]^{ \kappa, \kappa}}
\nc{\Inc}{{\compactN^{\uparrow\bbN}}}
\nc{\powInc}[1]{{\big(\Inc\big)^{#1}}}
\nc{\powFin}[1]{{\big(\Fin\big)^{#1}}}
\nc{\powPN}[1]{{\big(\PN\big)^{#1}}}
\nc{\NcompactN}{{\compactN^\bbN}}
\nc{\Uarrow}{\smash{\big\uparrow}}
\nc{\LE}{\preccurlyeq}
\nc{\GE}{\succcurlyeq}
\nc{\op}{\operatorname}
\nc{\im}{\op{Im}}
\nc{\Span}{\op{span}}
\nc{\maxfin}{\op{maxfin}}
\nc{\ran}{\op{range}}
\nc{\iso}{\cong}
\nc{\Madd}{{\M}^\star}
\nc{\cI}{\mathcal{I}}
\nc{\cJ}{\mathcal{J}}
\nc{\scrA}{\mathscr{A}}
\nc{\scrB}{\mathscr{B}}
\nc{\scrC}{\mathscr{C}}
\nc{\scrD}{\mathscr{D}}
\nc{\scrF}{\mathscr{F}}
\nc{\scrK}{\mathscr{K}}
\nc{\A}{\D\forall}
\nc{\B}{\mathrm{B}}
\nc{\cB}{\mathcal{B}}
\nc{\cZ}{\mathcal{Z}}
\nc{\bB}{\mathbf{B}}
\nc{\BS}{\mathbf{B}(\mathcal{S})}
\nc{\BF}{\mathbf{B}(\mathcal{F})}
\nc{\BU}{\mathbf{B}(\mathcal{U})}
\nc{\cSp}{\mathcal{S}^+}
\nc{\cFp}{\mathcal{F}^+}
\nc{\cUp}{\mathcal{U}^+}

\nc{\BG}{\B_\Ga}
\nc{\BL}{\B_\Lambda}
\nc{\BT}{\B_\Tau}
\nc{\BTstar}{\B_{\Tau^*}}
\nc{\BO}{\B_\Om}
\nc{\DO}{\cD_\Om}
\nc{\KO}{\cK_\Om}
\nc{\CG}{C_\Ga}
\nc{\CL}{C_\Lambda}
\nc{\CT}{C_\Tau}
\nc{\CTstar}{C_{\Tau^*}}
\nc{\CO}{C_\Om}
\nc{\COgp}{C_{\Om^{\grbl}}}
\nc{\CLgp}{C_{\Lambda^{\grbl}}}
\nc{\BOgp}{\B_{\Om}^{\grbl}}
\nc{\BLgp}{\B_{\Lambda^{\grbl}}}
\nc{\sfC}{\mathsf{C}}
\nc{\sfD}{\mathsf{D}}
\nc{\bD}{\mathbf{D}}
\nc{\Tau}{\mathrm{T}}
\nc{\cA}{\mathcal{A}}
\nc{\cK}{\mathcal{K}}
\nc{\cD}{\mathcal{D}}
\nc{\cE}{\mathcal{E}}
\nc{\cF}{\mathcal{F}}
\nc{\cS}{\mathcal{S}}
\nc{\cT}{\mathcal{T}}
\nc{\cG}{\mathcal{G}}
\nc{\cY}{\mathcal{Y}}
\nc{\J}{\mathcal{J}}
\nc{\cL}{\mathcal{L}}
\nc{\cM}{\mathcal{M}}
\nc{\cN}{\mathcal{N}}
\nc{\cH}{\mathcal{H}}

\nc{\scF}{\mathscr{F}}
\nc{\scH}{\mathscr{H}}

\nc{\Op}{\mathrm{O}}
\nc{\rmA}{\mathrm{A}}
\nc{\rmF}{\mathrm{F}}
\nc{\rmB}{\mathrm{B}}
\nc{\rmD}{\mathrm{D}}
\nc{\rmP}{\mathrm{P}}
\nc{\cC}{\mathcal{C}}
\nc{\cP}{\mathcal{P}}
\nc{\bbQ}{\mathbb{Q}}
\nc{\bbR}{\mathbb{R}}
\nc{\cU}{\mathcal{U}}
\nc{\cX}{\mathcal{X}}
\nc{\cQ}{\mathcal{Q}}
\nc{\Un}{\bigcup}
\nc{\cV}{\mathcal{V}}
\nc{\cR}{\mathcal{R}}
\nc{\tcR}{\tilde{\mathcal{R}}}
\nc{\cW}{\mathcal{W}}
\nc{\Z}{{\mathbb Z}}
\nc{\Impl}{\Rightarrow}
\long\def\forget#1\forgotten{\marginpar{\textcolor{green}{Forgetting...}}}
\nc{\ft}{\mathfrak{t}}
\nc{\fb}{\mathfrak{b}}
\nc{\fc}{\mathfrak{c}}
\nc{\fd}{\mathfrak{d}}
\nc{\fg}{\mathfrak{g}}
\nc{\oo}{\infty}
\nc{\fr}{\mathfrak{r}}
\nc{\fk}{\mathfrak{k}}
\nc{\bidi}{\mathfrak{bidi}}
\nc{\fu}{\mathfrak{u}}
\nc{\fh}{\mathfrak{h}}
\nc{\fp}{\mathfrak{p}}
\nc{\fj}{\mathfrak{j}}
\nc{\fs}{\mathfrak{s}}
\nc{\w}{\omega}
\nc{\x}{\times}
\nc{\Iff}{\Leftrightarrow}
\newcommand\comp{^{\text{\tt c}}}
\nc{\nin}{\notin}
\nc{\cat}{\hat{\ }}
\nc{\sub}{\subseteq}
\nc{\spst}{\supseteq}
\nc{\sm}{\setminus}
\nc{\as}{\subseteq^*}
\nc{\les}{\le^*}
\nc{\leinf}{\le^{\infty}}
\nc{\leS}{\le_S}
\nc{\leF}{\le_F}
\nc{\leU}{\le_U}
\nc{\rest}{\restriction}

\nc{\la}{\langle}
\nc{\ra}{\rangle}
\nc{\E}{\exists}
\nc{\dom}{\op{dom}}
\nc{\cov}{\op{cov}}
\nc{\add}{\op{add}}
\nc{\addmen}{\add(\Men{})}
\nc{\cof}{\op{cof}}
\nc{\cf}{\op{cf}}
\nc{\non}{\op{non}}
\nc{\unif}{\op{non}}
\nc{\COV}{\op{COV}}
\nc{\ADD}{\op{ADD}}
\nc{\COF}{\op{COF}}
\nc{\NON}{\op{NON}}
\nc{\impl}{\to}
\nc{\Lp}{\mathcal{L_\p}}
\nc{\Wlog}{without loss of generality}
\newtheorem{thm}{Theorem}[section]
\nc{\bthm}{\begin{thm}} \nc{\ethm}{\end{thm}}
\newtheorem{need}[thm]{Need}
\nc{\bneed}{\begin{need}\color{dg}} \nc{\eneed}{\end{need}}
\newtheorem{prop}[thm]{Proposition}
\nc{\bprp}{\begin{prop}} \nc{\eprp}{\end{prop}}
\newtheorem{fact}[thm]{Fact}
\nc{\bfct}{\begin{fact}} \nc{\efct}{\end{fact}}
\newtheorem{prob}[thm]{Problem}
\nc{\bprb}{\begin{prob}} \nc{\eprb}{\end{prob}}
\newtheorem{lem}[thm]{Lemma}
\nc{\blem}{\begin{lem}} \nc{\elem}{\end{lem}}
\newtheorem{app}[thm]{Application}
\nc{\bapp}{\begin{app}} \nc{\eapp}{\end{app}}
\newtheorem{claim}[thm]{Claim}
\nc{\bclm}{\begin{claim}} \nc{\eclm}{\end{claim}}
\newtheorem{cor}[thm]{Corollary}
\nc{\bcor}{\begin{cor}} \nc{\ecor}{\end{cor}}
\newtheorem{conj}[thm]{Conjecture}
\nc{\bcnj}{\begin{conj}} \nc{\ecnj}{\end{conj}}
\theoremstyle{definition}
\newtheorem{defn}[thm]{Definition}
\nc{\bdfn}{\begin{defn}} \nc{\edfn}{\end{defn}}
\newtheorem{obs}[thm]{Observation}
\nc{\bobs}{\begin{obs}} \nc{\eobs}{\end{obs}}
\theoremstyle{remark}
\newtheorem{rem}[thm]{Remark}
\nc{\brem}{\begin{rem}} \nc{\erem}{\end{rem}}
\newtheorem{cnv}[thm]{Convention}
\nc{\bcnv}{\begin{cnv}} \nc{\ecnv}{\end{cnv}}
\newtheorem{exam}[thm]{Example}
\nc{\bexm}{\begin{exam}} \nc{\eexm}{\end{exam}}
\nc{\bpf}{\begin{proof}} \nc{\epf}{\end{proof}
}
\nc{\be}{\begin{enumerate}}
\nc{\ee}{\end{enumerate}}
\nc{\bi}{\begin{itemize}}
\nc{\bimy}{\my{\begin{itemize}}
\nc{\eimy}{\end{itemize}}}
\nc{\itm}{\item}
\nc{\ei}{\end{itemize}}
\nc{\Subsection}[1]{\goodbreak\subsection*{#1}}

\nc{\sone}{\mathsf{S}_1}
\nc{\sfin}{\mathsf{S}_\mathrm{fin}}
\nc{\ufin}{\mathsf{U}_\mathrm{fin}}
\nc{\Split}{\mathsf{Split}}

\nc{\gone}{\mathsf{G}_1}    
\nc{\tgfin}{{\mathsf{G}}^*_\mathrm{fin}}
\nc{\gfin}{\mathsf{G}_\mathrm{fin}}

\nc{\men}[1]{\sfin(\Op(#1),\Op(#1))}
\nc{\sch}{\ufin(\cO,\Omega)}
\nc{\rothb}{\text{Rothberger}}
\nc{\pmen}{\sfin(\Omega,\Omega)}
\nc{\Rothb}{\sone(\Op,\Op)}

\nc{\prothb}{\sone(\Omega,\Omega)}
\nc{\tU}{{\tilde{U}}}
\nc{\tF}{{\tilde{F}}}
\nc{\tY}{{\tilde{Y}}}
\nc{\tX}{{\tsup[1]{X}}}
\nc{\dtX}{{\tsup[2]{X}}}
\nc{\dt}[1]{{\tsup[2]{#1}}}
\nc{\td}{{\tilde{d}}}
\nc{\tz}{{\tilde{z}}}
\nc{\cfd}{\cf(\fd)}

\nc{\msep}{\sfin(\cD,\cD)}
\nc{\rsep}{\sone(\cD,\cD)}
\nc{\cft}{\sfin(\Omega_{\mathbf{0}},\Omega_{\mathbf{0}})}
\nc{\scft}{\sone(\Omega_{\mathbf{0}},\Omega_{\mathbf{0}})}

\nc{\Umen}{U\text{-Menger}}
\nc{\hur}{\ufin(\cO,\Gamma)}
\nc{\tUmen}{\tU\text{-Menger}}

\nc{\Men}{\text{Menger}}
\nc{\Sch}{\text{Scheepers}}

\nc{\aspst}{\prescript{*}{}{\spst}\ }
\nc{\eqs}{=^*}
\nc{\ctblOm}{\Omega_{\mathrm{ctbl}}}
\nc{\GNga}{{\smallbinom{\Om}{\Ga}}}
\nc{\ctblga}{\smallbinom{\ctblOm}{\Ga}}

\nc{\nadd}{\cN_{\mathrm{add}}}
\nc{\ball}{\mathrm{B}}
\nc{\cOunif}{\cO^{\textrm{unif}}}

\nc{\sep}{
\vspace{2cm}
\noindent
\begin{minipage}{\textwidth}
	\textcolor{red}{\rule{\textwidth}{1pt}}
\end{minipage}
}
\nc{\FS}{\op{FS}}
\nc{\sums}{\op{SS}}
\nc{\SG}{\op{SG}}
\nc{\tSG}{\SG_\odot}
\nc{\G}{\op{G}}
\nc{\FSG}{\op{FSG}}
\nc{\FP}{\op{FP}}
\nc{\nonNadd}{\non(\nadd)}
\nc{\borga}{\Ga_\mathrm{Bor}}
\nc{\borw}{\Omega_\mathrm{Bor}}
\nc{\pick}{x}
\nc{\gen}{y}
\nc{\nullzind}{\sone(\{\Op_n^{\mathsf{unif}}\}_{n\in\bbN},\Ga)}
\nc{\nullzindf}[1]{\sone(\{\Op_{#1}^{\mathsf{unif}}\}_{n\in\bbN},\Ga)}

\definecolor{dg}{RGB}{42,101,24}

\nc{\myb}[1]{{\color{blue}{#1}}}
\nc{\mydg}[1]{{\color{dg}{#1}}}
\DeclareMathOperator{\eexists}{\exists}
\DeclareMathOperator{\fforall}{\forall}
\nc{\Exists}[1]{\bigl(\eexists #1\bigr)}
\nc{\Forall}[1]{\bigl(\fforall #1\bigr)}
\nc{\End}[1]{\bigl(#1\bigr)}
\nc{\dmo}[2]{\DeclareMathOperator{#1}{#2}}
\dmo{\Asc}{Asc}
\nc{\plusmin}{\wedge}

\nc{\cBsub}{{\cB^{\mbox{\tiny $\sub$}}}}

\nc{\Alice}{{\textsc{Alice}{}}}
\nc{\Bob}{{\textsc{Bob}}}

\nc{\bfO}{\mathbf{0}}

\nc{\proba}[1]{F}

\nc{\opn}[1]{\Op}
\nc{\sfinoo}[1]{\sfin(\Op,\Op)}
\nc{\om}[1]{\Om}

\nc{\restrict}{\upharpoonright}
\nc{\sgom}{\mathsf{S}_1(\Gamma,\Omega)}

\nc{\NNbar}{\overline{\w}^{\,\w}}
\nc{\lekap}{\le^\kappa_\kappa}
\nc{\rothkap}{[\kappa]^\kappa}
\nc{\Finkap}{\kappa^{<\kappa}}
\nc{\gakap}{\gamma_\kappa}
\nc{\wkap}{\w_\kappa}
\nc{\finkap}{\mathsf{Q}_\kappa}
\nc{\fbkap}{\fb_\kappa}
\nc{\fdkap}{\fd_\kappa}
\nc{\fpkap}{\fp_\kappa}
\nc{\ftkap}{\ft_\kappa}
\nc{\pkap}{P_\kappa}
\nc{\Gakap}{\Gamma_\kappa}
\nc{\U}[2]{U^{(#1)}_{#2}}
\nc{\subkap}{\sub^*_\kappa}
\nc{\Omkap}{\Omega_\kappa}
\nc{\Succ}{\mathrm{Succ}_\kappa}
\nc{\Opk}{\Opn_\kappa}
\nc{\Gak}{\Ga_\kappa}
\nc{\Omk}{\Om_\kappa}

\nc{\PK}{\op{P}(\kappa)}
\nc{\leskap}{\leq_\kappa^*}
\nc{\ZFC}{\mathsf{ZFC}}

\nc{\hurGEN}{\ufinGEN{\cO_\kappa}{\Gamma_\kappa}}
\nc{\sfinGEN}[2]{\mathsf{S}_{<\kappa}^\kappa(#1,#2)}
\nc{\ufinGEN}[2]{\mathsf{U}_{<\kappa}^\kappa(#1,#2)}
\nc{\soneGEN}[2]{\mathsf{S}_{1}^\kappa(#1,#2)}
\nc{\menGEN}{\sfinGEN{\cO_\kappa}{\cO_\kappa}}
\newcommand{\kk}{\kappa^\kappa}
\newcommand{\kkup}{{\kappa^{\uparrow \kappa}}}

\keywords{Hurewicz property, Menger property, Rothberger property, gamma property, selection principles, generalized Baire space, generalized Cantor space, cardinal characteristics, weakly compact cardinal}

\subjclass[2020]{Primary 54D20; Secondary 03E17, 03E55, 54A25, 54A35}

\author[A.~Amsalem]{Ayelet Amsalem}
\address{Ayelet Amsalem, Department of Mathematics, Bar-Ilan University, Ramat Gan 5290002, Israel}
\email{evyaye@gmail.com}

\author[A.~Jarden]{Adi Jarden}
\address{Adi Jarden, Department of Mathematics, Ariel University, Ariel 40700, Israel}
\email{jardena@ariel.ac.il}

\author[M.~Pawlikowski]{Micha\l{} Pawlikowski}\address{Micha\l{} Pawlikowski, Faculty of Technical Physics, Information
   Technology and Applied Mathematics,
   Lodz University of Technology, Aleje Politechniki 8,
       93-590 Ł\'od\'z}
\email{michal-pawlikowski4@wp.pl}

\author[P. Szewczak]{Piotr Szewczak}
\address{Institute of Mathematics, University of Warsaw, Banacha 2,
02–097 Warsaw, Poland}
\email{p.szewczak@wp.pl}
\urladdr{http://piotrszewczak.pl}

\author[B.~Tsaban]{Boaz Tsaban}
\address{Boaz Tsaban, Department of Mathematics, Bar-Ilan University, Ramat Gan 5290002, Israel}
\email{tsaban@math.biu.ac.il}
\urladdr{http://math.biu.ac.il/~tsaban}

\author[L. Zdomskyy]{Lyubomyr Zdomskyy}
\address{Institut f\"ur Diskrete Mathematik und Geometrie, Technische Universit\"at Wien, Wiedner Hauptstrasse 8-10/104, 1040 Wien, Austria.}
\email{lzdomsky@gmail.com}
\urladdr{https://dmg.tuwien.ac.at/zdomskyy/}

\thanks{
The research of Micha\l{} Pawlikowski and Piotr Szewczak was funded by the National Science Center, Poland Weave-UNISONO call in the Weave programme
Project: Set-theoretic aspects of topological selections 2021/03/Y/ST1/00122.
This research has been completed while Micha\l{} Pawlikowski was the Doctoral Candidate in the Interdisciplinary Doctoral School at the 
{\L}\'{o}d\'{z} University of Technology, Poland.
The research of Lyubomyr Zdomskyy
was funded in whole by the Austrian Science Fund (FWF) [10.55776/I5930 and 10.55776/PAT5730424].
}

\title[]{Combinatorial covering properties in an uncountable setting: canonical examples}

\begin{document}
\maketitle

\begin{abstract}
We provide examples of spaces satisfying generalized combinatorial covering properties such as the Hurewicz, Menger, and $\gamma$-properties in an uncountable setting.
Our approach is motivated by canonical constructions from the classical countable case, including the examples of Bartoszyński and Shelah separating the Hurewicz property from $\sigma$-compactness, the examples of Tsaban and Zdomskyy separating the Hurewicz and Menger properties, and Tsaban's construction of a nontrivial set of reals with the $\gamma$-property.
We focus on the genuinely nontrivial aspects of these higher-cardinal generalizations, uncovering several open problems whose nature appears substantially different from that of their countable counterparts.
\end{abstract}

\section*{User manual}

Combinatorial covering properties as the Hurewicz, Menger, and $\gamma$-properties~\cite{Menger2002, Hurewicz1926berEV, GN1982} are among the central notions studied in modern set-theoretic topology.
The primary objects of investigation are \emph{sets of reals}, that is, topological spaces homeomorphic to subspaces of the Cantor cube $\Cantor$, which are of independent interest.
These properties have also found applications in several areas, including forcing~\cite{CRZ}, local properties of function spaces~\cite{Ark, Sak, Buk, GN1982}, products of topological spaces~\cite{SzeGO,AT}, Ramsey theory in algebra~\cite{Tsa, SzeColor}, and hyperspace theory~\cite{Kru}.

During numerous conferences devoted to topology and set theory, we were repeatedly asked whether these properties admit meaningful generalizations to an uncountable setting, where the Cantor cube $\Cantor$ is replaced by its natural counterpart $2^\kappa$ for an uncountable cardinal $\kappa$.
Formulating the corresponding definitions presents little difficulty, and the first substantial progress in this direction was made by Korch--Weiss~\cite{KoWe}.
However, the construction of genuinely nontrivial examples of spaces with generalized covering properties by combinatorial methods in $2^\kappa$ seems to have received comparatively little attention.
This is somewhat surprising, as decades of research in the classical countable setting have produced a rich collection of elegant, albeit often highly nontrivial, constructions and examples.
Often, generalizations from the countable to the uncountable setting are essentially devoid of new mathematical content.
This happens when the proofs require no genuinely new ideas and amount merely to replacing $\w$ by an uncountable cardinal $\kappa$ in classical arguments.
Generalized combinatorial covering properties are easy to formulate, and their combinatorial characterizations are frequently straightforward adaptations of the countable case.
Our goal was to avoid such trivialities.
Accordingly, whenever an argument is merely a routine modification of a known proof, we omit the details.
Instead, we focus on genuinely nontrivial examples.
Although the underlying intuition often originates from classical constructions, obtaining the desired results requires new ingredients and considerably more delicate arguments.
In several cases, such as the $\kappa$-Hurewicz Conjecture (see Section~\ref{sec:HC}), this approach leads to highly nontrivial and apparently difficult open questions.

\section{Introduction}

Let $\kappa$ be an infinite cardinal number.
A topological space $X$ is a \emph{$P_\kappa$-space}, if intersection of less than $\kappa$-many open subsets of $X$ is open.
By \emph{space} we mean an infinite $P_\kappa$-space which is Tychonoff.
Let $\cU$ be a cover of a space $X$.
The family $\cU$ is a \emph{$\kappa$-cover} of $X$ if $X \notin \cU$ and for each set $F \sub X$ with $\card{F}<\kappa$ there is a set $U \in \cU$ such that $F \sub U$.
The family $\cU$ is a \emph{$\gakap$-cover} of $X$ if it has size $\kappa$ and the sets 
$\set{U \in \cU}{x \notin U}$ have size less than $\kappa$ for all $x\in X$.
By $\Opn_\kappa, \Omega_\kappa,$ and  $\Ga_\kappa$ we denote the classes of all \emph{open covers of size $\kappa$}, all \emph{open $\kappa$-covers of size $\kappa$} and all \emph{open $\gakap$-covers} of spaces, respectively.
If $\cA$ is a class of covers of spaces, then by $\cA(X)$ denote the family of all covers of $X$ from the class $\cA$.
Let $\cA,\cB \in \{\Opn_\kappa, \Omega_\kappa, \Ga_\kappa\}$.
We define the following generalized combinatorial covering properties which a space $X$ can possess.
\mbox{}\\

\begin{labeling}{$\ufinGEN{\cA}{\cB}$:}
\setlength{\itemsep}{0.5em}
\item [$\soneGEN{\cA}{\cB}$:]  for each $\seq{\cU_\alpha}{\alpha < \kappa}\in\cA(X)^\kappa$, there are  $U_\alpha\in\cU_\alpha$ such that $\set{U_\alpha}{\alpha<\kappa}\in \cB(X)$, 
\item [$\sfinGEN{\cA}{\cB}$:] for each $\seq{\cU_\alpha}{\alpha<\kappa}\in\cA(X)^\kappa$, there are  $\cV_\alpha \in [\cU_\alpha]^{<\kappa}$ such that $\Un_{\alpha<\kappa}\cV_\alpha \in \cB(X)$,
\item [$\ufinGEN{\cA}{\cB}$:] for each $\seq{\cU_\alpha}{\alpha<\kappa}\in\cA(X)^\kappa$, there are $\cV_\alpha\in [\cU_\alpha]^{<\kappa}$ such that $\set{\Un\cV_\alpha}{\alpha<\kappa}\in \cB(X)$.	
\end{labeling}
\mbox{}\\

The properties $\mathsf{U}_{<\w}^\w(\Opn,\Gamma)$, $\mathsf{S}_{<\w}^\w(\Opn,\Opn)$,  $\mathsf{S}_1^\w(\Om,\Ga)$ and $\mathsf{S}_1^\w(\Opn,\Opn)$, are classical celebrated Hurewicz~\cite{Hurewicz1926berEV}, Menger~\cite{Menger2002}, the $\gamma$-property~\cite{GN1982}, and Rothberger~\cite{Rothberger1938}, respectively.
The above properties were introduced by Korch and Weiss~\cite{KoWe} as generalizations of properties from the Scheepers diagram~\cite{coc1} to the case of sequences of covers of uncountable length.

Now we explain what we mean by "natural counterpart of $\Cantor$ in an uncountable setting".
Let $\kappa$ be an infinite regular cardinal and $X$ be a discrete topological space.
A \emph{$\kappa$-product topology} on $\prod_{\alpha < \kappa} X_\alpha$ is the one generated by the sets $\prod_{\alpha \in \kappa}U_\alpha$ where $U_\alpha\sub X_\alpha$ for all $\alpha \in \kappa$ and $U_\alpha= X$ for all but less than $\kappa$ many $\alpha$'s.
The main stream of investigations of these properties is devoted to spaces which are homeomorphic with subspaces of the \emph{$\kappa$-Cantor space} $2^\kappa$, endowed with the $\kappa$-product topology.
By a \emph{set} with a topological property we mean a space homeomorphic with a subspace of $2^\kappa$ with this property.
Identifying subsets of $\kappa$ with their characteristic functions, elements of $2^\kappa$, we introduce a topology on $\PK$. 
We identify each element $a \in \rothkap$ with the increasing enumeration of its elements, an element of the \emph{$\kappa$-Baire space} $\kappa^\kappa$ considered with the $\kappa$-product topology.
The topologies inherited from $2^\kappa$ and $\kappa^\kappa$ coincide on $\rothkap$.
Both the $\kappa$-Cantor space and $\kappa$-Baire space are Tychonoff $P_\kappa$-spaces that are zero-dimensional.
If $\kappa^{<\kappa} = \kappa$, then $\kappa$-Cantor space and $\kappa$-Baire space have bases of size $\kappa$. 
From now on we assume that $\kappa^{<\kappa} = \kappa$, and that $\kappa$ is uncountable and \emph{strongly inaccessible}, i.e., it is regular and $2^\alpha<\kappa$ for all $\alpha<\kappa$.
\mbox{}\\

One of the main streams of investigations in the countable case were motivated, by the following problems (and their solutions)~\cite{BS01, MR2421163, BT, FM, TsabanTheBook, OrTs}.

\bi
\item Menger's Conjecture: Is there in $\ZFC$ a Menger set of reals which is not $\sigma$-compact?
\item Hurewicz's Conjecture: Is there in $\ZFC$ a Hurewicz set of reals which is not $\sigma$-compact?
\item Hurewicz's Problem: Is there in $\ZFC$ a Menger set of reals which is not Hurewicz?
\item $\gamma$-set Problem: Is there (using very weak set-theoretic assumptions) a nontrivial $\gamma$-set of reals?
\ei
We address counterparts of the above problems in an uncountable setting.
A space is \emph{$\kappa$-compact} if every open cover of this space has a subcover of size less than
$\kappa$.
A space is $K_\kappa$ if it is a union of $\kappa$ many $\kappa$-compact spaces.
A space is $\kappa$-Hurewicz if it satisfies $\ufinGEN{\Opn_\kappa}{\Gakap}$, \emph{$\kappa$-Menger} if it satisfies $\sfinGEN{\Opn_\kappa}{\Opn_\kappa}$, and $\gamma_\kappa$ if any $\kappa$-cover has a $\gamma_\kappa$-subcover.
A cardinal number $\kappa$ is \emph{weakly compact} if 
$2^\kappa$ is $\kappa$-compact~\cite[8.23. Theorem]{MR396267}.
A weakly compact cardinal is strongly inaccessible.
The following diagram presents counterparts of the mentioned problems from the countable case.
In the case of dashed arrows we assume that $\kappa$ is weakly compact and in others we assume that  $\kappa$ is only strongly inaccessible.

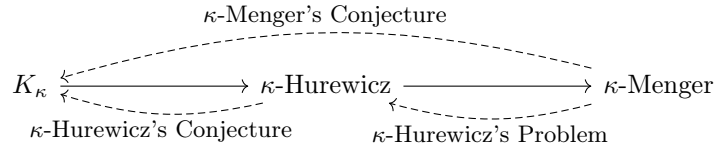
\begin{figure}[H]
\begin{tikzcd}[column sep=2.5cm]
K_\kappa
\arrow{r}
&
\arrow[bend right=-15, dashed]{dlu}{\text{\small $\kappa$-Hurewicz's Conjecture}}
\kappa\text{-Hurewicz}
\arrow{r}
&
\arrow[bend right=-15, dashed]{dlu}{\text{\small $\kappa$-Hurewicz's Problem}}
\arrow[bend left=-15,swap, dashed]{ulld}{\text{\small $\kappa$-Menger's Conjecture}}
\kappa\text{-Menger}
\end{tikzcd}
\end{figure}

Section~\ref{sec:GM} is devoted to the Galvin--Miller Lemma~\cite[Lemma~1.2]{MR738943} which is one of the fundamental tools used in constructions, and therefore it has been separated from other parts.
We also expect that it can find applications considering products of spaces --- as it was in the countable setting~\cite{unbddtower, ST,MR3449260}.
In Section~\ref{sec:HC} we focus on an examples due to Bartoszy\'{n}ski--Shelah~\cite{BS01} with a proof provided by Tsaban~\cite{TsabanTheBook}.
These examples show that in $\ZFC$, there are nontrivial Hurewicz sets of reals with no homeomorphic copies of the Cantor set inside.
Consequences of this part are the most surprising.
In contrast to the countable case, a structure of spaces from these adapted examples does not provide automatically that they are not $K_\kappa$.
Moreover, there are some indications that consistently such spaces can be $K_\kappa$ and also that the $\kappa$-Hurewicz Conjecture may be consistent with $\ZFC +\kappa$ is weakly compact.

Section~\ref{sec:HP} concerns the $\kappa$-Hurewicz Problem. The construction presented there is based on ideas from the final version of the Tsaban--Zdomskyy solution of the classical Hurewicz Problem, together with material from Amsalem's M.Sc. thesis, supervised by Jarden and Tsaban~\cite{MSc}.
An example provided in this section is also a solution of $\kappa$-Menger Conjecture.
Here we have another difference with the countable case, where a counterexample in the countable case for Menger's Conjecture is not on this level of difficulty~\cite{BT}.
Section~\ref{sec:gamma} is a construction of $\gamma_\kappa$-space using a counterpart of set-theoretic assumptions used by Tsaban in the countable case~\cite{OrTs, OST}.
We finish the paper with a list of questions and possible future directions.

\section{Galvin--Miller Lemma}\label{sec:GM}
One of the standard tools used in constructions in classical case is a technical Lemma due to Galvin--Miller~\cite[Lemma~1.2]{MR738943}.
The importance of this lemma suggests that it deserves for a separate section.

\bdfn
A function $f \in \kk$ is
\bi
\item \emph{strictly increasing} if for all $\alpha<\beta<\kappa$ 
we have $f(\alpha)<f(\beta)$,
\item \emph{continuous} if for any 
limit ordinal $\alpha<\kappa$ we have 
$f(\alpha) = \sup_{\beta<\alpha}f(\beta)$.
\item \emph{normal} if it is strictly increasing and continuous.
\ei
\edfn

\brem Observe that if $f \in \kk$ is normal then
\[
\Un_{\alpha<\kappa} [f(\alpha), f(\alpha+1)) = [f(0), \kappa).
\]
\erem

Define $\finkap:=[\kappa]^{<\kappa}$.

\blem[{Galvin--Miller~\cite[Lemma~1.2]{MR738943}}]\label{lem:Galvin_Miller_gen}
Let $\cU \in \Omega_\kappa(\finkap)$ be a family of open subsets of $\PK$. Then
there are a normal function $a \in \kk$ and a sequence $\seq{U_\alpha}{\alpha<\kappa}$ of pairwise distinct sets from $\cU$ 
such that for every $x \in \Pof(\kappa)$ and $\alpha<\kappa$:
\begin{equation}\label{eq:GM}
\text{If } 
x \cap [a(\alpha), a(\alpha+1)) = \emptyset,\text{ then }x \in U_\alpha.
\end{equation}
\elem

In our proof we need the following rather easy fact.

\bfct\label{lem:fin_set_containment} 
Let $X$ be a space with $\card{X}\geq \kappa$ and $\cU$ be an $\wkap$-cover of $X$.
Then for each set $A \sub X$ with $\card{A}<\kappa$, there are $\kappa$ many sets $U \in \cU$ such that $A \sub \cU$.
\efct 

\bpf[Proof of Lemma~\ref{lem:Galvin_Miller_gen}]
Let $a(0) := 0$ and $U_0\in\cU$ be a set such that $\Pof([0,a(0)))=\{\emptyset\}\sub U_0$.
Fix an ordinal number $\alpha<
\kappa$ and assume that $a(\alpha)<\kappa$ and sets $U_\beta\in \cU$ have been defined for all $\beta<\alpha$ such that for any limit ordinal $\alpha'\leq\alpha$ we have $a(\alpha')=\sup_{\beta<\alpha'}a(\beta)$.

Define $a(\alpha+1)$ and $U_\alpha\in\cU$ as follows.
Since $\kappa$ is strongly inaccessible, we have $\card{\Pof([0, a(\alpha)))}<\kappa$.
By Lemma~\ref{lem:fin_set_containment}, each subset of $\finkap$ with cardinality less than $\kappa$, is contained in $\kappa$ many 
elements of $\cU$, thus there is a set 
\[
U_\alpha \in \cU \sm \sset{U_\beta}{\beta < \alpha}\quad\text{with}\quad \Pof([0, a(\alpha))) \sub U_\alpha.
\]
Fix $s \in \Pof([0, a(\alpha)))$.
Since $s \in U_\alpha$ and the set $U_\alpha$ is open in $\PK$, there is $a(\alpha)<\alpha_s<\kappa$ such that 
\[
B_s:=\sset{s\cup t}{t\in \Pof([\alpha_s,\kappa))}\sub U_\alpha,
\]
the set $B_s$ is an open basic neighborhood of $s$ in $\PK$.
Let 
\[
a(\alpha +1):=\sup_{s\in \Pof([0,a(\alpha)))}\alpha_s.
\]
Since $\kappa$ is regular, we have $a(\alpha+1)<\kappa$.

Fix $x\in \PK$ such that $x\cap[a(\alpha), a(\alpha+1))=\emptyset$.
Then for $s:=x\cap[0,a(\alpha))$, we have $x\in B_s\sub U_\alpha$.
\epf

The Galvin--Miller Lemma has served as a sharp tool in many constructions related to products of spaces in the countable case~\cite{unbddtower, ST,MR3449260}. 
We expect that the above modification will also play an important role when considering products of spaces in the uncountable setting.

\section{$\kappa$-Hurewicz Conjecture}\label{sec:HC}

Hurewicz's Conjecture asserts that, in the class of metrizable spaces, $\sigma$-compactness is equivalent to the Hurewicz property.
This Conjecture was refuted by Just--Miller--Scheeprs--Szeptycki~\cite{coc1} (using a dichotomic proof) and then by Bartoszy\'{n}ski--Shelah~\cite{BS01} (providing a $\ZFC$ uniform counterexample).
The counterexample of Bartoszy{\'n}ski--Shelah is a starting point for our canonical examples.
We need the following notions.
For functions $a, b \in \roth$ we write $a \les b$, if the set 
$\sset{n \in \w}{a(n)> b(n)}$ is finite.
A set $A\sub\roth$ is \emph{bounded}, if there is a function $b\in\NN$ such that $a\les b$ for all $a\in A$.
Let $\fb$ be the minimal cardinality of a subset of $\roth$ which is \emph{unbounded}, i.e., it is not bounded, 
A set $X = \sset{x_\alpha}{\alpha<\fb} \sub \roth$ is a \emph{$\fb$-scale} if $X$ is unbounded  and $x_\alpha \les x_\beta$ for all $\alpha<\beta<\fb$.
A $\fb$-scale exists in $\ZFC$.
Let $\Fin:=[\w]^{<\w}$.
A set of the form $X \cup \Fin$, where $X$ is a $\fb$-scale is called a \emph{$\fb$-scale set}.
Bartoszyński--Shelah~\cite{BS01} proved that each $\fb$-scale set is Hurewicz but is not $\sigma$-compact.

The first item on our to do list is to extend above notions for combinatorics in $\rothkap$ (easy), verifying whether a counterpart of a $\fb$-scale set in $\PK$ is $\kappa$-Hurewicz (achievable after proper adjustments) and showing that each such a generalized $\fb$-scale set is not $K_\kappa$ (there are some indications that it is impossible to show this outright in ZFC).

For functions $a, b \in \rothkap$, we write $a \leskap b$ if the set $\sset{\alpha \in \kappa}{a(\alpha)> b(\alpha)}$ is of size less than $\kappa$.
We write $a \le^\kappa_\kappa b$ if the set $\sset{\alpha \in \kappa}{a(\alpha) \le b(\alpha)}$ is of size $\kappa$.
A set $A\sub\rothkap$ is \emph{bounded} if there is a function $b\in\rothkap$ such that $a\leskap b$ for all $a\in A$.
Let $\fbkap$ be the minimal cardinality of subset of $\rothkap$ which is \emph{unbounded}, i.e., it is not bounded.
The value of $\fbkap$ does not change if we consider its counterpart in $\kappa^\kappa$ instead of $\rothkap$.
A set $A\sub\rothkap$ is \emph{dominating} if for each function $x\in\rothkap$, there is a function $a\in A$ such that $x\leskap a$.
Let $\fdkap$ be the minimal cardinality of a dominating set in $\rothkap$.
The study of $\fbkap$ and $\fdkap$ has been initiated by Cummings and Shelah~\cite{MR1355135}
who investigated global behavior of these cardinal functions.

A set $X = \sset{x_\alpha}{\alpha<\fbkap} \sub \rothkap$ is a \emph{$\fbkap$-scale} if $X$ is unbounded  and $x_\alpha \leskap x_\beta$ for all $\alpha< \beta<\fbkap$.
It can be shown that $\fbkap$-scales exist in $\ZFC$. 
A combinatorial structure of these sets justifies that $\fbkap$ is a regular cardinal and $\kappa<\fbkap$.
A set of the form $X \cup \finkap$, where $X$ is a $\fbkap$-scale is called \emph{$\fbkap$-scale set}. 

A zero-dimensional metrizable space $X$ is Hurewicz if and only if $X$ is Lindel\"of and each continuous image of $X$ into $\roth$ is bounded~\cite{MR1271477}.
Basically, the same approach allows to get the following characterization. 
Since a proof has a one-to-one correspondence to the classical one~\cite{MR1271477}, it is omitted.
A space is $\kappa$-Lindel\"of if any cover of this space has a subcover of size at most $\kappa$. Let $P$ be a topological property such that there is a subspace of $\PK$ which does not have this property. 
Let $\non(P)$ be the minimal cardinality of a subspace of $\PK$ without the property $P$.

\bprp\label{prp:H_char}
Let $X$ be a zero-dimensional $P_\kappa$-space.
The space $X$ is $\kappa$-Hurewicz if and only if $X$ is $\kappa$-Lindel\"of and each 
continuous image of $X$ into $\rothkap$ is $\leskap$-bounded.
\eprp

\bcor
$\non(\kappa\text{-Hurewicz})=\fbkap$.
\ecor

In fact we can obtain even stronger result than in the above corollary, which will be useful in our considerations.

\blem\label{lem:spaces_smaller_than_bkap_are_S1(G,G)}
$\non(\soneGEN{\Gamma_\kappa}{\Gamma_\kappa}=\fbkap$.
\elem

\bpf
Assume that $X$ is $\kappa$-Lindel\"of, zero-dimensional $P_\kappa$-space.
Let $\seq{\cU_\alpha}{\alpha<\kappa}$ be a sequence of open $\gakap$-covers of $X$.
Since $X$ is $\kappa$-Lindel\"{o}ff, we may assume that $\cU_\alpha=\sset{\U{\alpha}{\xi}}{\xi<\kappa}$.
Let $\varphi \colon X \to \kk $ be a map such that $\varphi(x)(\alpha)<\kappa$ is the minimal ordinal such that 
\[
x \in \bigcap_{\varphi(x)(\alpha)\leq \xi<\kappa} \U{\alpha}{\xi}.
\]
Since $\card{X}<\fbkap$, the set $\varphi[X]$ is $\leskap$-bounded by a function $b\in\kk$.
Then $\sset{\U{\alpha}{b(\alpha)}}{\alpha < \kappa}$ is a $\gakap$-cover of $X$. 
\epf

\brem
In the above proof, for any function $c\in\kk$ with $b\leskap c$, the family $\sset{\U{\alpha}{c(\alpha)}}{\alpha<\kappa}$ is a $\gakap$-cover of $X$.
\erem

Another useful Corollary from Proposition~\ref{prp:H_char} is the following one.

\bcor \label{cor:union_of_kappa_many_Kappa_Hur_is_bounded}
Let $\sset{X_\xi}{\xi<\kappa}$ be a family of $\kappa$-Hurewicz spaces. 
Let $X = \Un_{\xi<\kappa} X_\xi$. If $X$ is a zero-dimensional $P_\kappa$-space, then  it is $\kappa$-Hurewicz. 
\ecor 

The following proof is motivated by a proof due to Tsaban~\cite[Theorem 2.12]{TsabanTheBook}.

\bthm\label{thm:bkscale}
Let $X$ be a $\fbkap$-scale.
Then $X \cup \finkap$ is $\kappa$-Hurewicz.
\ethm

\bpf
Let $X = \sset{x_\alpha}{\alpha<\fbkap}$ be a $\fbkap$-scale and $\seq{\cU_\alpha}{\alpha<\kappa}$ be a sequence of open covers of $X \cup \finkap$, which are open also in $\PK$ that are of size $\kappa$.
For each $\alpha<\kappa$, applying Lemma~\ref{lem:Galvin_Miller_gen} to the family $\cU_\alpha$, we get a function $a_\alpha$ and a sequence $\seq{U^{(\alpha)}_\beta}{\beta<\kappa}$ of sets from $\cU_\alpha$ such that
for every $x \in \Pof(\kappa)$ and $\alpha<\kappa$:
\[
\text{If } 
x \cap [a_\alpha(\beta), a_\alpha (\beta+1)) = \emptyset,\text{ then }x \in U^{(\alpha)}_\beta.
\]
For a function $c\in\kk$ such that $c(\xi) := a_\xi(\card{\xi}^+)$ for all $\xi<\kappa$, there is $\alpha<\kappa$ such that $x_\alpha\not\leskap c$.
Then a set 
\[
I := \sset{\xi<\kappa}{a_\xi(\card{\xi}^+) < x_\alpha(\xi)}
\]
has size $\kappa$.

By Lemma~\ref{lem:spaces_smaller_than_bkap_are_S1(G,G)}, the set $\card{\sset{x_\beta}{\beta<\alpha}}<\fbkap$ satisfies $\soneGEN{\Gamma_\kappa, \Gamma_\kappa}$.
Thus there is a function $b\in \kappa^I$ such that 
\[
\sset{\U{\xi}{b(\xi)}}{\xi \in I} \in \Gakap\left(\set{x_\beta}{\beta<\alpha}\right).
\]
Define 
\[
\cV_\xi :=\begin{cases}
\emptyset,\text{ if }\xi \in\kappa\sm I,\\
\sset{\U{\xi}{\beta}}{\beta \le \card{\xi}^+} 
\cup \{\U{\xi}{b(\xi)}\}, \text{ if }\xi \in I.
\end{cases}
\]

Fix $x \in X \cup \finkap$.
Assume that $x \in \finkap$.
Since the functions $a_\xi$ are strictly increasing for all $\xi<\kappa$, for each large enough $\xi \in I$ we have $x \cap [a_\xi(\xi), a_\xi(\xi+1))=\emptyset$. 
Then $x \in \U{\xi}{\xi} \sub \Un\cV_\xi$ for all
large enough $\xi \in I$.
Assume that $x=x_\beta$ for some $\beta<\alpha$.
Then $x_\beta \in \U{\xi}{b(\xi)} \sub\Un \cV_\xi$ for all large enough 
$\xi \in I$.
Now assume that $x=x_\beta$ for some $\beta$ with $\alpha\leq \beta<\kappa$.
Then for all but less than $\kappa$ many $\xi \in I$ we have 
\[
x_\beta(\xi) > x_\alpha(\xi) > a_\xi(\card{\xi}^+)>\card{\xi}^+,
\]
and thus 
\[
\card{x_\beta \cap \bigl[0, a_\xi(\card{\xi}^+)\bigr)} < \card{\xi}^+.
\]
Fix any $\xi$ witnessing this property. 
Since 
\[
[0, a_\xi(\card{\xi}^+)) = \Un_{\delta < \card{\xi}^+} [a_\xi(\delta), a_\xi(\delta+1))
\]
is a union of $\card{\xi}^+$ intervals, there is $\delta < \card{\xi}^+$ such that 
\[
x_\beta \cap [a_\xi(\delta), a_\xi(\delta+1)) = \emptyset.
\]
It follows that $x \in \Un \cV_\xi$ for all but less than $\kappa$ many $\xi \in I$.

We conclude that 
\[
\smallmedset{\Un\cV_\xi}{\xi<\kappa}\in\Gakap(X\cup\finkap).\qedhere
\]
\epf
Let $\lambda$ be an infinite cardinal number.
A set $X \sub \PK$ is 
\emph{$\lambda$-concentrated on a set $Q\sub X$} if for each open set $U\sub \PK$ containing $Q$ we have 
$\card{X \sm U}<\lambda$.

\bfct\label{lem:image_of_compact_is_compact}\label{lem:compact_image_bounded}
Let $X$ be a $\kappa$-compact space and
$f \colon X \to \kappa$ be a continuous function, where in $\kappa$ we consider the order topology.
Then $f[X]$ is $\kappa$-compact, and thus there is $\alpha<\kappa$ with $f[X] \sub \alpha$.
\efct

\blem 
Assume that $\kappa$ is a weakly compact cardinal.
Let $X$ be a $\fbkap$-scale. 
Then $X \cup \finkap$ is $\fbkap$-concentrated on $\finkap$.
\elem 
\bpf 
Let $X = \sset{x_\alpha}{\alpha<\fbkap}$ be a 
$\fbkap$-scale and $U \sub \PK$ be an 
open set containing $\finkap$.
Then the set $\PK \sm U \sub \rothkap$ 
is closed, and thus it is $\kappa$-compact. 
Fix $\alpha<\kappa$.
Define an 
evaluation map $e_\alpha \colon \rothkap \to \kappa$ such that $e_\alpha(a) =: a(\alpha)$.
For each $\beta<\kappa$ the set
\[
e_\alpha^{-1}[[0, \beta)) = 
\smallmedset{a \in \rothkap}{a(\alpha) \in [0,\beta)} = 
\Un_{\xi \in [0,\beta)} \set{a \in \rothkap}{a(\alpha) = \xi}
\]
is a union of less that $\kappa$ clopen subset of $\rothkap$, and thus it is clopen. 
It follows that the map $e_\alpha$ is 
continuous.

By Fact~\ref{lem:compact_image_bounded} for each
$\alpha<\kappa$, there exists $\beta_\alpha<\kappa$ such that 
$e_\alpha[\PK \sm U] \sub \beta_\alpha$.
Let $b\in\rothkap$ be a function such that $b(\alpha):= \beta_\alpha+1$. Then $a \leskap b$
for all $a \in \PK \sm U$.
Since the set $\sset{x_\alpha}{\alpha<\kappa}$ is $\leskap$-unbounded, there is $\alpha < \kappa$ such that $x_\alpha \not\leskap b$.
Then 
\[
X \sm U = X \cap (\PK \sm U) \sub 
\set{x_\beta}{\beta < \fbkap, x_\beta \leskap b} \sub 
\set{x_\beta}{\beta < \alpha}.\qedhere
\]
\epf 

\brem
In the countable case, one can show that a $\fb$-scale set in $\PN$ is $\fb$-concentrated on $\Fin$, and thus it does not contain an uncountable compact set.
Consequently, each $\fb$-scale set is not $\sigma$-compact.
In contrast to closed uncountable subsets of $\Cantor$, it is consistent with $\ZFC + \kappa$ is strongly inaccessible that there
is a closed subset $F \sub 2^\kappa$ with $\card{F}>\kappa $ such that there is no copy of $2^\kappa$ inside $F$~\cite{MR3235820}.
\erem

\section{$\kappa$-Hurewicz Problem}\label{sec:HP}

The original Hurewicz Problem asks whether in $\ZFC$, there is a Menger set of reals which is not Hurewicz~\cite{Hurewicz1927}.
It was first solved by Chaber--Pol~\cite{chaber2023remarkfremlinmillertheoremconcerning}(using a dychotomic proof) and then by Tsaban--Zdomskyy~\cite{MR2421163} (a~$\ZFC$ uniform construction). 
We briefly describe the latter construction here.
A set $A\sub\roth$ is \emph{dominating}, if for each function $x\in \roth$, there is a function $a\in A$ with $x\les a$.
Let $\fd$ be the minimal cardinality of a dominating set in $\roth$.
Let $\sset{d_\alpha}{\alpha<\fd}$ be a dominating set in $\roth$.
For $a,b\in\roth$, we write $a\leinf b$ if the set $\sset{n}{a(n)\leq b(n)}$ is infinite.
For $x \in \PN$ by $x\comp$ we denote the set $\w \sm x$. By $\ici$ we denote the family of all $x \in \roth$ such that $x\comp \in \roth$.
In $\ZFC$, there is a set $X=\sset{x_\alpha}{\alpha<\fd}\sub \ici$ such that 
\[
\sset{d_\beta}{\beta< \alpha}\leinf x_\alpha\quad\text{and}\quad d_\alpha\leinf x_\alpha\comp
\]
for all $\alpha<\fd$.
Then the set $X\cup\Fin$ is Menger but its homeomorphic image 
\[
\sset{x\comp_\alpha}{\alpha<\fd}\cup\sset{x\comp}{x\in\Fin}
\]
in $\roth$ is unbounded, and thus the set $X\cup\Fin$ is not Hurewicz.

Adjusting the Hurewicz Problem to our context we formulate the following question.

\bprb[$\kappa$-Hurewicz Problem]
Assume that $\kappa$ is a weakly compact cardinal number.
Is there a $\kappa$-Menger subspace of $\PK$ which is not $\kappa$-Hurewicz?

\eprb

A set $A\sub\rothkap$ is \emph{dominating} if for each function $x\in\rothkap$, there is a function $a\in A$ such that $x\leskap a$.
Let $\fdkap$ be the minimal cardinality of a dominating set in $\rothkap$.
A zero-dimensional metrizable space $X$ is Menger if and only if $X$ is Lindel\"of and each continuous image of $X$ into $\roth$ is not dominating~\cite[Proposition 3]{MR1271477}.
Basically, the same approach allows to get the following characterization of $\kappa$-Menger spaces. 
Since a proof has a one-to-one correspondence to the classical one, it is omitted.

\bprp\label{prp:M_char}
Let $X$ be a zero-dimensional $P_\kappa$-space.
The space $X$ is $\kappa$-Menger if and only if $X$ is $\kappa$-Lindel\"of and each 
continuous image of $X$ into $\rothkap$ is not dominating.
\eprp

\bcor
$\non(\kappa\text{-Menger})=\fdkap$.
\ecor

Our starting point is a combination of Tsaban's streamlined version of the Tsaban--Zdomskyy construction~\cite[Theorem~3.9]{MR2421163} and the corresponding uncountable analogue developed in Amsalem's M.Sc. thesis~\cite{MSc}.
The resulting combinatorial construction turns out to be substantially more delicate than its countable counterpart.

\bthm\label{thm:men}
Assume that $\kappa$ is weakly compact.
There is a $\kappa$-Menger subspace of $\PK$ which is not $\kappa$-Hurewicz.
\ethm

Before we prove Theorem~\ref{thm:men} we need additional technical results.

By $\Succ$ we denote the set of all successor ordinals smaller than $\kappa$.
\blem\label{lem:normal_func}
For each $f \in \kk$ there exists a normal function $g \in \kk$ such that 
$f\upharpoonright \Succ \le g\upharpoonright \Succ$.
\elem 
\bpf 
We may assume that $f$ is strictly increasing.
Let $g(0) := f(0)$.
Fix an ordinal number $\alpha$ with $0<\alpha<\kappa$ and assume that the values $g(\beta)$ have been defined for all $\beta<\alpha$.
If $\alpha$ is a successor ordinal then define $g(\alpha) := f(\alpha)$.
If $\alpha$ is a limit ordinal, then let 
\[
g(\alpha) := \sup_{\beta<\alpha} g(\beta) = \sup_{\beta < \alpha} f(\beta).\qedhere
\]
\epf 

\bdfn
For $x \in \rothkap$ define
\[
\cl{x}:=x \cup \sset{\alpha<\kappa}{\sup(\alpha \cap x)=\alpha}.
\]
\edfn

\blem\cite{MR1940513}\label{lem:en_of_club}
Let $x \in \rothkap$. The set $x$ is a $\kappa$-club if and only if
increasing enumeration of elements of $x$ is a normal function.
\elem 
\brem 
Lemma~\ref{lem:en_of_club} shows that for all $x \in \rothkap$, the 
function $\cl{x}$ is normal.
\erem

\blem\label{lem:big_partition}
Let $Y \sub \rothkap$ be a set such that $\card{Y}<\fdkap$.
Then there is a function
$b \in \rothkap$ such that 
\[
\card{\smallmedset{\alpha < \kappa}{\card{y \cap [b(\alpha), b(\alpha+1))} \ge 5}}=\kappa
\]
for all $y\in Y$.
\elem 

\bpf 
We may assume that for each set $y\in Y$, after removing less than $\kappa$ many elements from $y$, the resulting set is also in $Y$.
Suppose contrary to our assertion that for each $b \in \rothkap$ there is $y \in Y$ such that
\[
\smallmedcard{y \cap [b(\alpha), b(\alpha+1))} \le 4
\]
for all $\alpha<\kappa$.
Let $f\in\kk$ be a function such that
\[
f(\alpha):=
\begin{cases}
0,&\text{ if }\alpha=0,\\
f(\beta)+4,&\text{ if }\alpha=\beta+1,\\
\alpha,&\text{ if }\alpha\text{ is a limit ordinal}.
\end{cases}
\]
for all $\alpha<\kappa$.

Fix $a \in \rothkap$.
By Lemma~\ref{lem:normal_func} there is a continuous function $g \in \rothkap$ such that $a \rest\Succ \le g \rest\Succ$.
In particular, we have
\[
a(\alpha) \le a(\alpha+1) \le g(\alpha+1)
\]
for all $\alpha<\kappa$.
Pick $y \in Y$ such that 
\[
y \cap [0, g(0)) = \emptyset\quad\text{ and }\quad\card{y \cap [g(\alpha), g(\alpha+1))}\le 4
\]
for all $\alpha<\kappa$.

We show by induction that $g(\alpha) \le y(f(\alpha+1))$ for all $\alpha<\kappa$:
Since $y \cap [0, g(0))=\emptyset$, we have $g(0)\leq y(0)<y(f(1))$.
Fix $\alpha<\kappa$ and assume that $g(\beta) \le y(f(\beta+1))$ for all $\beta<\alpha$.

Assume that $\alpha=\beta+1$ for some $\beta<\kappa$.
We have 
\[
y(f(\alpha+1))=y(f(\beta+2))=y(f(\beta+1)+4),
\]
and thus $y(f(\alpha+1))$ is the fifth smallest element in the set $y\cap [y(f(\beta+1),\kappa)$.
Since 
\[
g(\beta)\leq y(f(\beta+1))\quad\text{ and }\quad\card{y\cap [g(\beta),g(\beta+1))}\leq 4,
\]
we have 
$g(\alpha)=g(\beta+1)\leq y(f(\alpha+1))$.

Assume that $\alpha$ is a limit ordinal.
Then
\[
g(\alpha) = \sup_{\beta<\alpha}g(\beta)
\le \sup_{\beta<\alpha}y(f(\beta+1))
\le y\bigl(\sup_{\beta<\alpha} f(\beta+1)\bigr)=y(\alpha)\leq
y(f(\alpha+1)).
\]

It follows that for each function $a\in\rothkap$ there is a function $y\in Y$ such that $a$ is dominated by 
\[
y'=\sset{y(f(\alpha+1))}{\alpha<\kappa}\in\rothkap.
\]
Since $\card{Y}<\fdkap$, we get a contradiction.
\epf 
The following is an immediate consequence of the previous lemma.

\bcor\label{cor:big_partition_open_intervals}
Assume that $\kappa$ is a strongly inaccessible cardinal.
Let $Y \sub \rothkap$ be a set such that $\card{Y}<\fdkap$.
Then there is 
$b \in \rothkap$ such that for each $y \in Y$, we have 
\[
\smallmedcard{\smallmedset{\alpha < \kappa}{\card{y \cap (b(\alpha), b(\alpha+1))} \ge 4}}=\kappa.
\]
\ecor 

Let $y \in \rothkap$ be a function such that $y\sub\Succ$.
Define 
$\tilde{y} \in \rothkap$ as follows:

\[
\tilde{y}(\alpha) := 
\begin{cases}
    y(0) + 1, \text{ if } \alpha = 0;\\
    y(\tilde{y}(\beta))+1,\text{ if } \alpha = \beta + 1; \\
    \max\{y(\alpha), \sup_{\beta<\alpha} \tilde{y}(\beta)\}, \text{ if }\alpha \text{ is a limit ordinal.}\\
\end{cases}
\]

By the definition $\tilde{y}$ is indeed an element of $\rothkap$.
We have $y\leq \tilde{y}$.

For $x, a \in \rothkap$, where $a$ is continuous define 
\[
x/a := \sset{\alpha<\kappa}{x \cap [a(\alpha), a(\alpha+1)) \neq \emptyset},
\]
an element of $\rothkap$.

\blem\label{lem:relationship_between_a_and_coml_of_a}
For each $a \in \rothkap$, we have:
\be 
\itm If $a(\alpha) = \overline{a}(\alpha')$ for some $\alpha,\alpha'<\kappa$, then $a(\alpha+1) = \overline{a}(\alpha'+1)$; 
\itm $\overline{a}(\alpha+1) \in a$ for all $\alpha<\kappa$.
\ee 
\elem 

\bpf
(1) Fix $\alpha, \alpha'<\kappa$ such that $a(\alpha) = \overline{a}(\alpha')$.
Let $\gamma \in (a(\alpha), a(\alpha+1))$ be a limit ordinal.
We have $\gamma\cap a = (a(\alpha) \cap a) \cup \{a(\alpha)\}$.
Then $\sup(\gamma \cap a) = a(\alpha) \neq \gamma$ and $\gamma \notin \overline{a}$.
It follows that $\overline{a}(\alpha'+1) > \gamma$ for all $\gamma \in (a(\alpha), a(\alpha+1))$.
We have $\overline{a}(\alpha'+1) \geq a(\alpha+1)$.
Since $a(\alpha+1) \in a$, we get $\overline{a}(\alpha'+1) \le a(\alpha+1)$, and thus $\overline{a}(\alpha'+1) = a(\alpha+1)$.

(2) Fix $\alpha<\kappa$.
If $\overline{a}(\alpha) \in a$, then by~(1) we have $\overline{a}(\alpha+1) \in a$.
Assume that $\overline{a}(\alpha) \notin a$.
Then $\overline{a}(\alpha) = \sup_{\beta<\beta_0} a(\beta) < a(\beta_0)$ for some $\beta_0<\kappa$. 
Let $\gamma \in (\overline{a}(\alpha), a(\beta_0))$ be a limit ordinal.
Then $\gamma \cap a = \overline{a}(\alpha) \cap a$.
Consequently $\sup(\gamma \cap a) = \overline{a}(\alpha) \neq \gamma$, and thus $\gamma \notin \overline{a}$. It follows  that $\overline{a}(\alpha+1) = a(\beta_0) \in a$.
\epf

\blem\label{lem:existence_of_a_function_jumping_over_a_family_kappa_times}
Let $Y \sub \rothkap$ be a set with $\card{Y}<\fdkap$ and $a \in \rothkap$.
Then there exists $I \in \rothkap$ such that
\[
Y\lekap \Un_{\alpha \in I} [a(\alpha), a(\alpha+1)).
\]
\elem 

\bpf 
Since the set $[\Succ]^\kappa$ is $\leskap$-dominating in $\rothkap$, we may assume that $Y\sub[\Succ]^\kappa$.

Let $a \in \rothkap$.
Define
\[
\tilde{Y}/\overline{a} \coloneqq \sset{\tilde{y}/\overline{a}}{y \in Y}.
\]
Let $b \in \rothkap$ be a function from Corollary~\ref{cor:big_partition_open_intervals}, applied to the set $\tilde{Y}/\overline{a}$.  
Fix $y \in Y$.
Then the set
\[
J:=\smallmedset{\alpha<\kappa}{\smallmedcard{\tilde{y}/\overline{a} \cap \bigl(b(\alpha), b(\alpha+1)\bigr)} \geq 4}
\]
has size $\kappa$.
Note that $b(\alpha)+2 < b(\alpha+1)$ for all $\alpha \in J$.
Firstly, we prove that for
\[
d' \coloneqq \Un_{\alpha \in b} [\overline{a}(\alpha), \overline{a}(\alpha+2)).
\]
we have $Y\le^\kappa_\kappa d'$:

Fix $\alpha \in J$ and four consecutive points $\gamma_1<\gamma_2<\gamma_3<\gamma_4$ in 
\[
\tilde{y}/\overline{a} \cap \bigl(b(\alpha), b(\alpha+1)\bigr).
\] 
We have
\[
b(\alpha)<\gamma_1<\gamma_2<\gamma_3<\gamma_4<b(\alpha+1).
\]
Define 
\[
\beta:= \sup(\tilde{y} \cap [\overline{a}(\gamma_2), \overline{a}(\gamma_2+1)))\quad\text{and}\quad
\beta':= \min\sset{\gamma<\kappa}{\tilde{y}(\gamma)\geq \beta}.
\]
Note that if $\beta \in \tilde{y}$, then $y(\beta')=\beta$.
Since $\gamma_1, \gamma_2, \gamma_3, \gamma_4$ are proximate in 
$\tilde{y}/\overline{a} \cap \bigl(b(\alpha), b(\alpha+1)\bigr)$, we have
\[
\tilde{y}(\beta')\in [\overline{a}(\gamma_2)\cup \overline{a}(\gamma_2+1))\cup[\overline{a}(\gamma_3), \overline{a}(\gamma_3+1))
\]
and
\[
\tilde{y}(\beta'+1)\in[\overline{a}(\gamma_3), \overline{a}(\gamma_3+1))\cup[\overline{a}(\gamma_4), \overline{a}(\gamma_4+1)).
\]
It follows that 
\[
y(\beta) \le y(\tilde{y}(\beta')) \le \tilde{y}(\beta'+1) <\overline{a}(\gamma_4+1) \le \overline{a}(b(\alpha+1)).
\]
We have also
\[
\overline{a}(b(\alpha)+2) \le \overline{a}(\gamma_2) \le \beta \le \overline{a}(\gamma_2+1)<\overline{a}(b(\alpha+1)).
\]
By the definition of $d'$, we have
\[
d' \cap [\overline{a}(b(\alpha)+2), \overline{a}(b(\alpha+1))) = \emptyset.
\]
The inequalities $\overline{a}(b(\alpha)+2)\leq\beta \le d'(\beta)$ imply that $\overline{a}(b(\alpha+1)) \le d'(\beta)$.
We conclude that  
\[
y(\beta) < \overline{a}(b(\alpha+1)) \le d'(\beta).
\]
It follows that $Y \le^\kappa_\kappa d'$.

Now fix $\alpha \in J$.
By Lemma~\ref{lem:relationship_between_a_and_coml_of_a}, we have $a\cap [\overline{a}(\alpha),\overline{a}(\alpha+2))\neq\emptyset$.
Then there exists $\beta<\kappa$ such that $a(\beta) \in \{\overline{a}(\alpha), \overline{a}(\alpha+1)\}$.
By Lemma~\ref{lem:relationship_between_a_and_coml_of_a} we have 
$a(\beta+1) \in \{\overline{a}(\alpha+1), \overline{a}(\alpha+2)\}$, and thus
\[
[a(\beta), a(\beta+1)) \sub [\overline{a}(\alpha), \overline{a}(\alpha+2)) \sub d'.
\] 
Then the set 
\[
I:=\sset{\alpha<\kappa}{[a(\alpha),a(\alpha+1))\sub d'}
\]
has size $\kappa$.
Since
\[
d := \Un_{\alpha \in I} [a(\alpha), a(\alpha+1))\sub d'
\]
we have $d'\le d$, and thus $Y \le^\kappa_\kappa d$.
\epf 

\blem\label{lem:jumping_over_arbitrary_function_by_x_and_x^c} 
Let $d \in \rothkap$.
Then there exists $a \in \rothkap$, such that for each 
$I \in\icikap$
letting 
\[x \coloneqq \Un_{\alpha \in I} [a(\alpha), a(\alpha+1))\] we have 
$d\lekap x$ and $d \lekap x\comp$.
\elem
\bpf 
Since the set $[\Succ]^\kappa$ is $\lekap$-dominating in $\rothkap$, we may assume that $d \in [\Succ]^\kappa$.
Set $a = \tilde{d}$, fix $I\in\icikap$ and
let $\alpha \in I\comp$.
We have $x \cap [\tilde{d}(\alpha), \tilde{d}(\alpha+1)) = \emptyset$.
Since $\tilde{d}(\alpha) \le x(\tilde{d}(\alpha))$, we have $\tilde{d}(\alpha+1) \le x(\tilde{d}(\alpha))$, and thus
\[
d(\tilde{d}(\alpha))
<
d(\tilde{d}(\alpha))+1
=
\tilde{d}(\alpha+1)
\le
x(\tilde{d}(\alpha)).
\]
It follows that $d \le^\kappa_\kappa x$.

Now let $\alpha \in I$. Since $x\comp(a(\alpha)) \geq a(\alpha)$ and 
$[a(\alpha), a(\alpha+1)) \cap x\comp = \emptyset$. Thus 
\[
d(a(\alpha)) = d(\tilde{d}(\alpha)) \le \tilde{d}(\alpha+1) = a(\alpha+1) \le x\comp(a(\alpha)).
\]
It follows that $d \lekap x\comp$.
\epf 

\blem\label{lem:constructing_special_dominating_set}
Let $Y\sub\rothkap$ be a set such that $\card{Y}<\fdkap$, and $d\in \rothkap$.
Then there is $x\in \icikap$ such that 
\[
Y\lekap x\quad \text{ and }\quad d\lekap x\comp.
\]
\elem
\bpf 
Let $a\in\rothkap$ be a function from Lemma~\ref{lem:jumping_over_arbitrary_function_by_x_and_x^c}, applied to $d$.
By Lemma~\ref{lem:existence_of_a_function_jumping_over_a_family_kappa_times},
there is $I \in \rothkap$ such that 
\[Y \lekap \Un_{\alpha \in I} [a(\alpha), a(\alpha+1)).\]
Eventually, by skipping to a subset of $I$, we may assume that $I\in\icikap$.
Let 
\[x \coloneqq \Un_{\alpha \in I} [a(\alpha), a(\alpha+1)).\]
Then $Y \lekap x$ and by Lemma~\ref{lem:jumping_over_arbitrary_function_by_x_and_x^c} 
we have $d \lekap x\comp$.
\epf

Let $\lambda$ be an infinite cardinal number. 
A set $X \sub \rothkap$ is \emph{$\lambda$-unbounded} if 
$\card{X} \ge \lambda$ and for each $g \in \rothkap$ we have $\card{\sset{x \in X}{x \leskap g}}<\lambda$.

\blem\label{lem:mu_unbounded_and_concetrated}
Assume that $\kappa$ is a weakly compact cardinal.
Let $\lambda$ be an infinite cardinal number and $X \sub \rothkap$ be 
$\lambda$-unbounded set.
Then $X \cup \finkap$ is $\lambda$-concentrated on $\finkap$.
\elem 

\bpf 
Let $U$ be an open set in $\PK$ containing $\finkap$.
Then $\PK\sm U\sub\rothkap$ is a closed subset of 
$\PK$.
By the assumption, $\PK$ is $\kappa$-compact, and thus $\PK\sm U$ is a $\kappa$-compact, too.
There exists $g \in \rothkap$ such that $\PK\sm U \le g$.
It follows that $X \sm U\sub\sset{x \in X}{x \le g}$.
Since $X$ is $\lambda$-unbounded, the latter set has size smaller than $\lambda$.
\epf 

\bcor\label{clm:fdkap_concentrated_is_Menger}
Assume that $\kappa$ is a weakly compact cardinal.
If $X \sub \rothkap$ is a $\fdkap$-unbounded set, then the set $X\cup\finkap$ is $\fdkap$-concentrated on $\finkap$, and thus $X\cup\finkap$ is $\kappa$-Menger.
\ecor 

\bpf[{Proof of Theorem~\ref{thm:men}}]
Let $\sset{d_\alpha}{\alpha<\fdkap}$ be a dominating set in $\rothkap$.
By Lemma~\ref{lem:constructing_special_dominating_set}, there exists a set
\[
X \coloneqq \sset{x_\alpha}{\alpha<\fdkap}\sub\icikap
\]
such that 
\[
\sset{d_\beta}{\beta\le\alpha} \lekap x_\alpha\quad \text{ and }\quad d_\alpha \lekap x_\alpha^\mathsf{c}
\]
for all $\alpha<\kappa$.

Fix $\alpha<\fdkap$.
Since $\sset{d_\beta}{\beta\le\alpha} \lekap x_\alpha$, we have 
$\sset{x\in X}{x \leskap d_\alpha} \sub \sset{x_\beta}{\beta<\alpha}$, and thus $\card{\sset{x\in X}{x \leskap d_\alpha}}<\fdkap$.
Since the set $\sset{d_\alpha}{\alpha<\fdkap}$ is dominating, the set $X$ is $\fdkap$-unbounded.
By Lemma~\ref{lem:mu_unbounded_and_concetrated} and Corollary~\ref{clm:fdkap_concentrated_is_Menger}, the set $X\cup\finkap$ is $\kappa$-Menger.

A function $\varphi\colon \PK\to\PK$ such that $\varphi(x):=x\comp$ is a homeomorphism and $\varphi[X\cup\finkap]\sub\rothkap$.
Since the set $\sset{d_\alpha}{\alpha<\fdkap}$ is dominating and $d_\alpha \lekap x_\alpha\comp$, the set $\sset{x_\alpha\comp}{\alpha<\fdkap}$ is unbounded.
Since $\sset{x_{\alpha}^\mathsf{c}}{\alpha<\fdkap} \sub \varphi[X \cup \finkap]$, the set $\varphi[X \cup \finkap] \sub \rothkap$ is unbounded, too.
It follows that a continuous image of $X\cup\finkap$ into $\rothkap$ is unbounded.
By Proposition~\ref{prp:H_char}, the set $X\cup\finkap$ is not $\kappa$-Hurewicz.
\epf

\section{$\kappa$-Menger Conjecture}\label{sec:MC}

The original Menger Conjecture asserts that in the class of metrizable spaces, $\sigma$-compactness is equivalent to the Menger property.
This conjecture was refuted firstly, by Fremlin--Miller~\cite{FM} (using a dichotomic proof) and then by Bartoszy{\'n}ski--Tsaban~\cite{BT} (a uniform $\ZFC$ construction).
A counterexample provided by Bartoszy{\'n}ski--Tsaban has the following form. Let $\sset{d_\alpha}{\alpha<\fd}$ be a dominating set in $\roth$.
Then there is a set $X=\sset{x_\alpha}{\alpha<\fd}\sub\roth$ such that 
\[
\sset{d_\beta}{\beta<\alpha}\leinf x_\alpha
\]
for all $\alpha<\fd$.
Then the set $X\cup \Fin$ is Menger but it is $\fd$-concentrated on $\Fin$, and thus it does not contain an uncountable compact sets.
Consequently, it is not $\sigma$-compact.
Of course, we can formulate also an analogue of the Menger Conjecture.

\bcnj[$\kappa$-Menger Conjecture]
Assume that $\kappa$ is a weakly compact cardinal.
A subspace of $\PK$ is a $K_\kappa$-space if and only if it is  $\kappa$-Menger.
\ecnj

A counterpart of the construction due to Bartoszy{\'n}ski--Tsaban can be achieved easily by transfinite induction also in the uncountable case.
However, a resulting set $X\cup\finkap$ is $\fdkap$-concentrated on $\finkap$ (and thus $\kappa$-Menger), but as it was already pointed out in the case of the $\kappa$-Hurewicz Conjecture, it does not provide that it is not $K_\kappa$.
The main result from the previous section provides a counterexample for the $\kappa$-Menger Conjecture, and we are not aware of any direct (not mentioning the $\kappa$-Hurewicz property) way to construct a $\kappa$-Menger non $K_\kappa$-space in $\ZFC$.

\section{nontrivial space with the $\gakap$-property}\label{sec:gamma}

Recall a non-trivial example of a space with the $\gamma$-property 
$\sone(\Omega, \Ga)$, given by Orenshtein and Tsaban~\cite[Theorem 3.6.]{OrTs}.
For $a, b \in \PN$ we write $a \sub^\ast b$ if the set $a \sm b$ 
is finite.
A family $A \sub [\w]^\w$ is \emph{centered} if for each finite nonempty subfamily $B \sub A$, we have $\card{\bigcap \cB}=\w$.
A \emph{pseudointersection} of $A \sub \roth$ is an element $b \in \roth$ such that $b \sub^\ast a$ for all $a\in A$, denoted by $b\as A$.
Let $\fp$ be the minimal cardinality of a centered subfamily of $\roth$ with no pseudointersection.
Let $\lambda$ be an infinite ordinal.
A set $X=\sset{x_\alpha}{\alpha<\lambda}\sub\roth$ is a \emph{tower} if 
for all $\alpha < \beta < \lambda$ we have $x_\beta \sub^\ast x_\alpha$ and $X$
does not have a pseudointersection.
Let $\ft$ be the minimal cardinality of a tower.
By a classical result of Maliaris--Shelah~\cite{MR3402699}, we have $\fp=\ft$.
A set $X\sub \roth$ is an \emph{unbounded tower} if it is a tower and it is $\les$-unbounded.

By the result of Orenshtein--Tsaban~\cite{OrTs}, if $X$ is an unbounded tower of cardinality $\fp$, then $X \cup \Fin$ has the $\gamma$-property.
An unbounded tower of cardinality $\fp$ exists if and only if $\fp = \fb$~\cite[Lemma 3.3]{OrTs}.
Now we proceed with defining analogues of pseudointersection, centeredness and towers.

For sets $a$ and $b$ we write $a \subkap b$ if $\card{a \setminus b}<\kappa$.
A family $A \sub \rothkap$ is \emph{$\kappa$-centered}, if for every nonempty family 
$B \sub A$ with $\card{B}<\kappa$ we have
$\card{\bigcap B}=\kappa$.
A \emph{$\kappa$-pseudointersection} of $A \sub \rothkap$ is an element
$b \in \rothkap$ such that $b \subkap a$ for all $a \in A$; denote this fact by $b\subkap A$.
Let $\fpkap$ be the minimal cardinality of a $\kappa$-centered subset of $\rothkap$ with no $\kappa$-pseudointersection.

Let $\lambda$ be an ordinal number.
A set $X = \sset{x_\alpha}{\alpha<\lambda}\sub\rothkap$
is a \emph{$\kappa$-tower} if $x_\beta \subkap x_\alpha$ for all $\alpha<\beta<\lambda$, $X$ is $\kappa$-centered and $X$ does not have a $\kappa$-pseudointersection.
Let $\ftkap$ be the minimal cardinality of a $\kappa$-tower.
A set $X\sub \rothkap$ is an \emph{unbounded $\kappa$-tower} if it is a $\kappa$-tower and it is $\leskap$-unbounded.
As in a countable case, it is possible to show the following relations between $\fpkap$ and $\ftkap$.
\blem[{\cite[Theorem 2.9]{MR4057946}}]
\label{lem:p_smaller_t}
$\kappa^+\le\fp_\kappa \le \ftkap \le \fbkap \le \fdkap \le 2^\kappa$. 
\elem 

The following fact connects the $\gakap$-property with the generalized combinatorial covering properties.
\bthm[Korch--Weiss~\cite{KoWe}]
A space $X \sub \PK$ with $\card{X} \geq \kappa$ has the $\gakap$-property if and only if it has the property $\soneGEN{\Omega_\kappa} {\Gamma_\kappa}$.
\ethm
The main goal of this section is to prove the following result.

\bthm\label{thm:gamma}
For each unbounded $\kappa$-tower $X \sub \rothkap$ of size $\fp_\kappa$, the space $X \cup \finkap$ is a $\gakap$-space.
\ethm 

A zero-dimensional metrizable space has the $\gamma$-property if and only if any continuous and centered image of the space into $\roth$ has a pseudointersection.
Firstly, we will present a useful combinatorial characterization of the $\gakap$-property which is a counterpart of the mentioned above.

\bprp
Let $X \sub \PK$ be a space.
Then $X$ is a $\gakap$-space if and only if each continuous and $\kappa$-centered 
image of $X$ into $\rothkap$ has a $\kappa$-pseudointersection. 
\eprp 

\bcor 
$\non(\gakap)=\fpkap$.
\ecor 


\bthm\label{thm:clubs_int}
The intersection of fewer than $\kappa$ 
closed unbounded subsets of $\kappa$ is 
closed unbounded.
\ethm

\brem 
The assumption of $\kappa$-centeredness in the definition of a $\kappa$-tower is necessary to obtain Lemma~\ref{lem:p_smaller_t}.
Without it we would obtain that $\ftkap=\w$ which is shown by the following example: 
The sets 
\[
x_n := \sset{(\lambda+n, \lambda + \w)}
{\lambda<\kappa \text{ is a limit ordinal}}
\]
are unbounded in $\kappa$ for all $n<\w$.
However $\bigcap_{n < \w} x_n = \emptyset$, thus the family $\sset{x_n}{n \in \w}$ does not have a $\kappa$-pseudointersection.

Furthermore the same example shows that  we cannot weaken the assumption of Theorem~\ref{thm:clubs_int} as it is easy to see that the sets $x_n$ are not closed.
\erem

\blem\label{lem:closure_of_pseudoint}
Let $A\sub\rothkap$ be a family of $\kappa$-clubs and $x\in\rothkap$ be a $\kappa$-pseudointersection of $A$.
Then $\cl{x}$ is a $\kappa$-pseudointersection of $A$.
\elem 
\bpf 
Let $a \in A$.
Since $x$ is a $\kappa$-pseudointersection of $A$,
there exists $\xi<\kappa$ such that  
$x \cap (\xi, \kappa) \sub a$. It suffices to show 
that for each limit ordinal $\alpha \in (\xi, \kappa)$ with $\sup(x \cap \alpha)=\alpha$,
we have $\alpha \in a$. Let $\alpha \in (\xi, \kappa)$ be such an ordinal. 
Since $\sup(x \cap \alpha)=\alpha$, for each $\beta \in (\xi, \alpha)$, there exists 
\[
\gamma \in x \cap (\beta, \alpha) \sub a \cap (\beta, \alpha),
\]
and thus
$\sup(a \cap \alpha) = \alpha$. 
Since $a$ is a $\kappa$-club, we have $\alpha \in a$.
\epf 

\blem\label{lem:bounding_tower}
Let $\lambda\le\ftkap$ and $\sset{y_\alpha}{\alpha<\lambda} \sub \rothkap$. Then there exists a set 
$\sset{x_\alpha}{\alpha<\lambda}\sub\rothkap$ such that 
for all $\alpha<\lambda$, we have 
\be 
\itm $x_\alpha$ is a $\kappa$-club; 
\itm $y_\alpha \upharpoonright \Succ \le x_\alpha \upharpoonright \Succ$;
\itm If $\beta<\alpha$, then 
$x_\alpha \subkap x_\beta$. 
\ee  
\elem 

\bpf 
Proceed by recursion.
Suppose we have defined $\sset{x_\beta}{\beta<\alpha}$ for some $\alpha<\lambda$.

Assume that $\cof(\alpha)<\kappa$.
Let $\seq{\beta_{\xi}}{\xi < \cof(\alpha)}$ 
be a cofinal sequence in $\alpha$.
By Theorem~\ref{thm:clubs_int} we have that 
$x \coloneqq \bigcap_{\xi<\cof(\alpha)} x_{\beta_\xi}$ is a $\kappa$-club.
Observe that $x$ is a $\kappa$-pseudointersection of $\sset{x_\beta}{\beta<\alpha}$.
By Lemma~\ref{lem:normal_func} there exists 
$x' \in \rothkap$ such that $x'$ is a $\kappa$-club and $y_\alpha \upharpoonright \Succ \le x' \upharpoonright \Succ$.
Then
$x_\alpha: = x' \cap x$ is a $\kappa$-club, it is a $\kappa$-pseudointersection 
of $\sset{x_\beta}{\beta<\alpha}$ and
$y_\alpha \upharpoonright \Succ \le x_\alpha \upharpoonright \Succ$.

Assume that $\cof(\alpha) \ge \kappa$.
We have $x_{\beta'}\as x_\beta$ for all $\beta<\beta'<\alpha$.
Since $\alpha<\ftkap$, there exists
a $\kappa$-pseudointersection $x$ of 
$\sset{x_\beta}{\beta<\alpha}$.
Then $\cl{x}$ is 
a $\kappa$-club and by Lemma~\ref{lem:closure_of_pseudoint} it is a 
$\kappa$-pseudointersection of
$\sset{x_\beta}{\beta<\alpha}$.
By Lemma~\ref{lem:normal_func} there exists 
$x' \in \rothkap$ such that $x'$ is a $\kappa$-club and $y_\alpha \upharpoonright \Succ \le x' \upharpoonright \Succ$.
Then $x_\alpha \coloneqq \cl{x} \cap x'$ is a $\kappa$-club set, it is is a $\kappa$-pseudointersection 
of $\sset{x_\beta}{\beta<\alpha}$ and
$y_\alpha \upharpoonright \Succ \le x_\alpha \upharpoonright \Succ$.
\epf 

\blem\label{lem:pseudoint_and_boundedness}
Let $X \sub \rothkap$ and $a \in \rothkap$ be such that the set $\sset{x \upharpoonright a}{x \in X}$ is $\leskap$-unbounded in $\kappa^a$.
Then $X$ has no $\kappa$-pseudointersection.
\elem

\bpf 
Suppose there exists $y \in \rothkap$ with $y\subkap X$.
Let $z\in\rothkap$ be a function such that $z(\alpha) := y(\alpha + \alpha)$ for all $\alpha<\kappa$.
Fix $x \in X$.
Since $y \subkap x$, there is $\alpha_0<\kappa$ such that $y(\alpha_0 + \alpha) \geq x(\alpha)$ for all $\alpha$ with $\alpha_0<\alpha<\kappa$.
Then $z(\alpha) \geq y(\alpha_0 + \alpha) \geq x(\alpha)$ for all $\alpha>\alpha_0$.
It follows that $X\leskap z$.
In particular $(x \upharpoonright a)  \leskap (z\restriction a)$ for all $x \in X$, a contradiction.
\epf 

\blem\label{lem:ext_of_undbb_ft_tower}
An unbounded $\kappa$-tower of size $\ftkap$ exists if and only 
if $\ftkap = \fbkap$.
\elem 

\bpf 
($\Rightarrow)$
The existence of an unbounded $\kappa$-tower of size $\ftkap$ implies that $\fbkap \le \ftkap$. Apply Lemma~\ref{lem:p_smaller_t}.

($\Leftarrow)$
Let 
$\sset{y_\alpha}{\alpha<\fbkap}\sub\rothkap$ be a $\leskap$-unbounded set. Let $\sset{x_\alpha}{\alpha<\ftkap}$ be a set from Lemma~\ref{lem:bounding_tower}.
The set $\sset{x_\alpha \upharpoonright \Succ}{\alpha<\ftkap}$ is $\leskap$-unbounded.
By Lemma~\ref{lem:pseudoint_and_boundedness} the set $\sset{x_\alpha}{\alpha<\ftkap}$ has no $\kappa$-pseudointersection. Furthermore since $\ftkap=\fbkap$ and $\fbkap$ is regular, the set $\sset{x_\alpha}{\alpha<\ftkap}$ is $\kappa$-centered.
Thus it is an unbounded $\kappa$-tower of size $\ftkap$.
\epf 

\blem\label{lem:kappa_boundedness_char}
Let $Y \sub \rothkap$.
The following assertions are equivalent:
\be 
\item The set $Y$ is $\leskap$-bounded;
\item There is $s \in \rothkap$ such that for each $a \in Y$ there is $\alpha_0<\kappa$ with 
\[
a \cap [s(\alpha), s(\alpha+1)) \neq \emptyset
\]
for all $\alpha \in (\alpha_0, \kappa)$.
\ee
\elem

\bpf 
($\Rightarrow)$
Let $b\in\rothkap$ be such that $Y\leskap b$.
By recursion, define $s \in \rothkap$ in the following way.
Put $s(0) := b(0)$.
Fix $\alpha<\kappa$ and assume that $s(\beta)$ has been defined for all $\beta<\alpha$.
Let $s(\alpha) := b(\sup_{\beta<\alpha} s(\beta))+1$.

Fix $a \in Y$.
Since $s \in \rothkap$, there 
exists $\alpha_0<\kappa$ such that $a(s(\alpha)) \le b(s(\alpha))$ for all $\alpha>\alpha_0$. It follows that
\[
s(\alpha) \le a(s(\alpha)) \le b(s(\alpha)) < b(s(\alpha))+1 = s(\alpha+1), 
\]
and thus $a \cap [s(\alpha), s(\alpha+1)) \neq \emptyset$ for all $\alpha$ with $\alpha<\alpha_0<\kappa$.

($\Leftarrow)$
Let $p_\alpha:=\sset{s(\alpha+\xi+1)}{\xi<\kappa}$ for all $\alpha<\kappa$ and definie $P:=\sset{p_\alpha}{\alpha<\kappa}$.
Since $\card{P}=\kappa$, there is $b\in\rothkap$ such that $P\leskap b$.
Fix $a\in Y$ and find $\alpha_0<\kappa$ such that $a\cap[s(\alpha),s(\alpha+1))\neq\emptyset$ for all $\alpha$ with $\alpha_0\leq\alpha<\kappa$.
Fix $\xi<\kappa$. 
Define 
\[
\beta(\xi):=\min\sset{\beta}{a(\beta)\in [s(\alpha_0+\xi),s(\alpha_0+\xi+1))}
\]
and note that $\beta(\xi)\geq \xi$.
It follows that 
\[
a(\xi)\leq a(\beta(\xi))\leq s(\alpha_0+\xi+1)=p_{\alpha_0}(\xi).
\]
We conclude that $a\leq p_{\alpha_0}\leskap b$.
\epf

\blem\label{lem:gamma_tower_lemma} 
Let $X \sub \PK$ be a space such that $\finkap \sub X$ and $\card{X}<\fp_\kappa$.
Let $\cU$ be a family of open sets in $\PK$ such that 
$\cU \in \Omega_\kappa(X)$ and $Y\sub\rothkap$ be a $\leskap$-unbounded set.
Then there are a set $a \in Y$ and a family of sets 
$\sset{U_\alpha}{\alpha < \kappa} \sub \cU$ such that 
$\sset{U_\alpha}{\alpha < \kappa} \in \Gamma_\kappa(X)$ and for each element
$x \in \rothkap$ and all $\alpha<\kappa$ we have:
\[
\text{If } x \sm \alpha \sub a, \text{ then } x \in \bigcap_{\xi \geq \alpha} U_\xi.
\]
\elem 
\bpf 
Since $\card{X}<\fp_\kappa$, it is a $\gakap$-space. Let 
$\cV \in \Gamma_\kappa(X)$ be a subfamily of $\cU$. 
By Lemma~\ref{lem:Galvin_Miller_gen} there are a function $b \in \rothkap$
and a family $\sset{V_{\alpha}}{\alpha<\kappa}$ of distinct sets
such that for each element $x \in \PK$ and all $\alpha<\kappa$:
\[
\text{If }
x \cap [b(\alpha), b(\alpha+1)) = \emptyset,\text{ then } x \in V_\alpha.
\]
By Lemma~\ref{lem:kappa_boundedness_char} there is a set $a \in Y$ such that the set
\[
c \coloneqq \sset{\alpha}{a \cap [b(\alpha), b(\alpha+1)) = \emptyset}
\]
is of size $\kappa$.

Fix $\alpha<\kappa$.
Take $\xi$ with  $\alpha\leq\xi<\kappa$ and any 
$x \in \PK$ such that $x \sm \alpha \sub a$.
Then $\alpha \leq c(\alpha) \le c(\xi)$ and 
\[
x \cap [b(c(\xi)), b(c(\xi+1))) \sub 
a \cap [b(c(\xi)), b(c(\xi+1))) = \emptyset.
\]
It follows that $x \in V_{c(\xi)}$ for all $\xi$ with  $\alpha\leq\xi<\kappa$, i.e., we have $x \in \bigcap_{\xi \geq \alpha} V_{c(\xi)}$. 
Define $U_\alpha \coloneqq V_{c(\alpha)}$ for all $\alpha < \kappa$.
Then $\sset{U_\alpha}{\alpha < \kappa} \in \Gamma_\kappa(X)$ and for each element 
$x \in \rothkap$ and all $\alpha < \kappa$ 
\[
\text{If }x \sm \alpha \sub a,\text{ then }x \in \bigcap_{\xi \geq \alpha} V_{c(\xi)} = \bigcap_{\xi \geq \alpha} U_\xi.\qedhere
\]
\epf 

\blem\label{lem:existence_of_unbdd_tower}
An unbounded $\kappa$-tower of size $\fp_\kappa$ exists if and only if $\fp_\kappa = \fbkap$.
\elem 
\bpf 
($\Rightarrow$)
Assume that an unbounded $\kappa$-tower of size $\fpkap$ exists. Then $\fp_\kappa \ge \fbkap$. By Lemma~\ref{lem:p_smaller_t} we have 
$\fp_\kappa \le \fbkap$.

($\Leftarrow$)
Now assume that $\fp_\kappa = \fbkap$. Then
$\fpkap=\ftkap = \fbkap$ and by Lemma~\ref{lem:ext_of_undbb_ft_tower}
there exists an unbounded $\kappa$-tower of size $\ftkap$.
\epf

\bpf[Proof of Theorem~\ref{thm:gamma}]
Let $X=\sset{x_\alpha}{\alpha<\fbkap} \sub \rothkap$ be a $\kappa$-tower that is 
unbounded with respect to $\leskap$. 
Let $X' = X \cup \finkap$ and for ordinal numbers 
$\gamma < \fbkap$ we define 
$X_\gamma = \sset{x_\alpha}{\alpha < \gamma} \cup \finkap$. Let $\cU \in \Omega_\kappa(X')$. Fix an ordinal number 
$\gamma_0 < \fbkap$.
Now we proceed by induction on $\beta \in (0, \kappa)$. Suppose we defined $\gamma_\xi$ for all $\xi<\beta$. Since
$\beta < \kappa$, it follows that 
$\gamma_\beta' =\sup_{\xi<\beta}\gamma_\xi<\kappa$. By Lemma~\ref{lem:gamma_tower_lemma} there 
are an ordinal number $\gamma_\beta<\fbkap$ and a subfamily $\sset{U_{\beta, \alpha}}{\alpha<\kappa} \in \Gamma_\kappa(X_{\gamma_\beta'})$ of $\cU$ such that for each element $x \in \rothkap$ and
all $\alpha<\kappa$ we have the following:
\[
\text{If } x \sm \alpha \sub x_{\gamma_\beta} \text{ then } x \in 
\bigcap_{\xi \ge \alpha}U_{\beta, \xi}.
\]
Let $\gamma = \sup_{\xi < \kappa} \gamma_\xi<\fbkap$. There is a function $g \in \kkup$ such that 
$x_\gamma \sm g(\xi) \sub x_{\gamma_\xi}$
for all $\xi<\kappa$.
Fix an ordinal number $\alpha \in [\gamma, \fbkap)$. Since 
$x_\alpha \subkap x_\gamma$ we have 
\[
x_\alpha \sm g(\beta) \sub x_\gamma \sm g(\beta) \sub x_{\gamma_\beta}
\]
for all but less than $\kappa$ many ordinals $\beta$.
Thus $x_\alpha \in \bigcap_{\xi \geq g(\beta)} U_{\beta, \xi}$ for all but less than $\kappa$ many $\beta$.
Thus for any function $h \in \kkup$
such that $g \le_{\kappa} h$ we have 
$\sset{U_{\xi, h(\xi)}}{\xi<\kappa} \in \Gamma_\kappa(X' \sm X_\gamma)$.

For each $x \in X_\gamma$ and each $\delta<\kappa$ we define 
\[
f_x(\delta) \coloneqq \min \set{\beta<\kappa}{x \in \bigcap_{\xi \geq \beta}U_{\delta, \xi}}
\]
if the set is nonempty and $f_x(\delta) = 0$ otherwise.
Note that if $\nu<\gamma$, then $\nu<\gamma_\delta$ for some $\delta<\kappa$, and thus there exists $\beta<\kappa$ with $x_\nu\in \bigcap_{\xi\geq \beta} U_{\delta,\xi}$, because $\sset{U_{\delta,\xi}}{\xi<\kappa}\in\Gamma_\kappa(X_{\gamma_\delta})$.

Since $\card{X_\gamma}<\fbkap$, there
is $h \in \kkup$ such that $h$ is an upper bound of $\{g\} \cup \sset{f_x}{x \in X_\gamma}$ with respect to $\leskap$ and the sets $U_{\beta, h(\beta)}$ are distinct. Then $\sset{U_{\beta, h(\beta)}}{\beta<\kappa} \in \Gamma_\kappa(X_\gamma)$. Since 
$\sset{U_{\beta, h(\beta)}}{\beta<\kappa} \in \Gamma_\kappa(X' \sm X_\gamma)$, as well, it follows that 
$\sset{U_{\beta, h(\beta)}}{\beta<\kappa} \in \Gamma_\kappa(X')$.
\epf 

\section{Example of a non-trivial space with the $\soneGEN{\Ga_\kappa}{\Ga_\kappa}$ property}

First we recall an example of a non-trivial set of reals with the classical $\sone(\Ga,\Ga)$ property given by Tsaban--Miller~\cite{MR2653961}, where $\sone(\Ga,\Ga)$ is the property $\sone^\w(\Gamma_\w,\Gamma_\w)$.
If $X\sub\roth$ is an unbounded tower of cardinality $\fb$, then $X \cup \Fin$ satisfies
$\sone(\Ga, \Ga)$. 
An unbounded tower of cardinality $\fb$ exists for example, if $\fp=\fb$ or $\fb<\fd$~\cite{MR2653961}
\bthm 
Let $X$ be an unbounded $\kappa$-tower of size $\fbkap$. 
Then $X \cup \finkap$ satisfies $\soneGEN{\Gamma_\kappa}{\Gamma_\kappa}$.
\ethm 
\bpf 
Let $X=\sset{x_\alpha}{\alpha<\fbkap}$ and 
$X_\alpha := \sset{x_\beta}{\beta<\alpha} \cup \finkap$ for all  $\alpha<\fbkap$.
Let $\seq{\cU_\alpha}{\alpha<\kappa}$ be a sequence of $\gakap$-covers of $X \cup \finkap$ which are open in $\PK$.
Fix $\alpha<\kappa$.
Applying Lemma~\ref{lem:Galvin_Miller_gen} to the family $\cU_\alpha$, we get a function $a_\alpha \in \rothkap$ and a sequence 
$\seq{U^{(\alpha)}_\xi}{\xi<\kappa}$ of sets from $\cU_\alpha$ such that for every $x \in \PK$ and $\xi<\kappa$ 
\[
\text{If }x \cap [a_\alpha(\xi), a_\alpha (\xi+1))=\emptyset,\text{ then } 
x \in U_\xi^{(\alpha)}.
\]

By Lemma~\ref{lem:kappa_boundedness_char}, for each $\alpha<\kappa$, there is $\delta_\alpha<\fbkap$ such that the set
\[
I_\alpha' := \sset{\xi<\kappa}{x_{\delta_\alpha} \cap [a_\alpha(\xi), a_\alpha(\xi+1)) = \emptyset}
\]
has size $\kappa$.
Let 
$\delta := \sup_{\alpha<\kappa} \delta_\alpha$. Since $\fbkap$ is regular, it follows that $\delta<\fbkap$. Furthermore for each $\alpha<\kappa$ we have that 
\[
x_\delta \cap [a_\alpha(\xi), a_\alpha(\xi+1)) \sub x_{\delta_\alpha} \cap [a_\alpha(\xi), a_\alpha(\xi+1))
\]
for all but less than $\kappa$ many 
$\xi<\kappa$.
Thus for each $\alpha<\kappa$ and all but less than $\kappa$ many $\xi \in I_\alpha'$ we have 
\[
x_\delta \cap [a_\alpha(\xi), a_\alpha(\xi+1)) = \emptyset.
\]
It follows that the sets
\[
I_\alpha := \sset{\xi<\kappa}{x_\delta \cap [a_\alpha(\xi), a_\alpha(\xi+1)) = \emptyset}
\] 
have size $\kappa$ for all $\alpha<\kappa$.

Since $\card{X_\delta}<\fbkap$,
there is a function $f \in \rothkap$
such that for any function $g\in\rothkap$ with $f\le g$, we have 
\[
\smallmedset{U^{(\alpha)}_{g(\alpha)}}{\alpha<\kappa}\in \Gakap(X_\delta).
\]
Take a function $g\in \rothkap\cap\prod_{\alpha<\kappa} I_\alpha$ with $f\le g$.

Fix $x \in X\sm X_\delta$.
Since $x \subkap x_\delta$, there is $\tilde{\alpha}<\kappa$ such that
\[
x\cap [\tilde{\alpha},\kappa)\sub x_\delta\cap[\tilde{\alpha},\kappa).
\]
Fix $\alpha$ with $\tilde{\alpha}\leq \alpha<\kappa$.
Since $g(\alpha)\in I_\alpha$, we have
\[
x\cap [a_\alpha(g(\alpha)),a_\alpha(g(\alpha)+1))\sub x_\delta \cap [a_\alpha(g(\alpha)), a_\alpha(g(\alpha) + 1)) = \emptyset,
\]
and thus $x\in U^{(\alpha)}_{g(\alpha)}$.

We conclude that 
\[
\sset{U^{(\alpha)}_{g(\alpha)}}{\alpha<\kappa}\in\Gakap(X\cup\finkap).\qedhere
\]
\epf 
An unbounded $\kappa$-tower of size $\fbkap$ exists for example if 
$\ftkap=\fbkap$.

Using Lemma~\ref{lem:pseudoint_and_boundedness} we can repeat the known argument to obtain that if $\fbkap<\fdkap$ then there is an unbounded $\kappa$-tower of size $\fbkap$.

\blem 
Assume that $\fbkap < \fdkap$ .
Then there is an unbounded $\kappa$-tower of 
size $\fbkap$.
\elem 

\bpf 
Let $X' \coloneqq \sset{x'_\alpha}{\alpha<\fbkap}$ be a 
$\fbkap$-scale. 
Since $\card{X'}<\fdkap$, the set $X'$ is not $\leskap$-
dominating.
Then there is $g \in \rothkap$ such that $x \lekap g$ 
for all $x \in X'$.
Define $x_\alpha \coloneqq \sset{\xi<\kappa}{x'_\alpha(\xi) \le g(\xi)}$. We will show that $X \coloneqq \sset{x_\alpha}{\alpha<\fbkap}$ is an unbounded $\kappa$-tower. 
Observe that $x_\beta \subkap x_\alpha$ for all $\alpha < \beta$, where $\alpha, \beta<\kappa$. By Lemma~\ref{lem:pseudoint_and_boundedness} it remains to show that
$X$ is $\leskap$-unbounded. 

Suppose that $X$ is $\leskap$-bounded.
Let $f \in \rothkap$ be a $\leskap$-bound of $X$. Then 
\[
x_\alpha'(\xi) \le x'_\alpha(x_\alpha(\xi)) \le g(x_\alpha(\xi)) \le g(f(\xi))
\]
for all but less than $\kappa$ many $\xi<\kappa$. 
It follows that $X'$ is $\leskap$-bounded, a contradiction.
\epf

\section{Comments and open problems}

\bprb
Is the $\kappa$-Hurewicz Conjecture  consistent with $\ZFC+\kappa$ is weakly compact?
\eprb

\bprb
Is it consistent with $\ZFC+\kappa$ is weakly compact that there is a $\fb_\kappa$-scale set which is $K_\kappa$?
\eprb

\bprb
Is there an alternative construction providing a $\kappa$-Menger but not $K_\kappa$-space than the one provided for the solution of $\kappa$-Hurewicz Problem?
\eprb

\bprb
What is the appearance of the Scheepers diagram in an uncountable case? How about the separations of other properties (even consistently)?
\eprb

\bibliography{bibliography.bib}
	 \bibliographystyle{abbrv}
\end{document}